\documentclass[12pt]{article}
\usepackage{amssymb}
\usepackage{latexsym}
\usepackage{amsmath}
\usepackage{graphicx}

\newenvironment{proof}[1][Proof]{\textbf{#1.} }{\ \rule{0.5em}{0.5em}}

\newtheorem{thm}{Theorem}[section]

\newtheorem{lem}[thm]{Lemma}
\newtheorem{prop}[thm]{Proposition}
\newtheorem{defn}{Definition}[section]
\newtheorem{rem}{Remark}[section]
\numberwithin{equation}{section}

\def\lbegar{$$\left\{ \begin{array}{lll}}
\def\rendar{\end{array} \right.$$}
\def\rendarp{\end{array} \right..$$}
\def\rendarv{\end{array} \right.,$$}

\newcommand\xih{\left(\frac{X_i-x}{h}\right)}

\newcommand\barO{\bar{B}}

\newcommand\lpl{LP(l) }

\newcommand\kxi{K\left(\frac{X_i-x}{h}\right)}
\newcommand\hth{\hat{\theta}}
\newcommand\sigmav{\vec{\sigma}}
\newcommand\lbeta{\lfloor \beta \rfloor} 

\def\px{P_X}
\def\eps{\epsilon}
\def\begar{$$\begin{array}{lll}}
\def\endar{\end{array}$$}
\def\begarlab{\begin{equation} \begin{array}{lll} \label}
\def\endarlab{\end{array} \end{equation}}

\def\ds1{{\mathrm{1 \hspace{-2.6pt} I}}}

\def\dsE{\mathbb {E}}
\def\dsN{\mathbb {N}}
\def\dsP{\mathbb {P}}

\def\dsR{\mathbb {R}}

\def\dsV{\mathbb {V}}
\def\dsZ{\mathbb {Z}}
\def\calA{{\cal A}}
\def\calB{{\cal B}}

\def\calC{{\cal C}}

\def\calFXY{{\cal F}}

\def\calN{{\cal N}}
\def\calP{{\cal P}}

\def\calS{{\cal S}}

\def\calX{{\cal X}}

\def\calZ{{\cal Z}}

\def\bg1{{g_1}}

\def\Z1N{Z_1^N}

\def\X1N{X_1^{N}}

\def\lam{{\lambda}}

\def\l0{{\lambda_0'}}

\def\Var{{\dsV\text{ar}} \,}

\def\hf{\hat{f}_n}
\def\sign{\text{sign}}

\def\ve{\varepsilon}

\newcommand\kap{\kappa}

\newcommand\demi{\frac{1}{2}}
\newcommand\fbeta{\lfloor \beta \rfloor}
\newcommand\heta{\hat{\eta}_n}

\newcommand\calPS{\calP_\Sigma}
\newcommand\calPSI{\calP_{\Sigma,\infty}}

\newcommand\Pn{P^{\otimes n}}
\begin{document}

\sloppy

\title{\bf
Fast learning rates for plug-in classifiers under the margin
condition }
\author{Jean-Yves AUDIBERT$^{1}$ and Alexandre B. TSYBAKOV$^2$\\
$^1$Ecole Nationale des Ponts et Chauss\'{e}es, $^2$Universit\'{e}
Pierre et Marie Curie} \maketitle

\begin{abstract}
It has been recently shown that, under the margin (or low noise)
assumption, there exist classifiers attaining fast rates of
convergence of the excess Bayes risk, i.e., the rates faster than
$n^{-1/2}$. The works on this subject suggested the following two
conjectures: (i) the best achievable fast rate is of the order
$n^{-1}$, and (ii) the plug-in classifiers generally converge slower
than the classifiers based on empirical risk minimization. We show
that both conjectures are not correct. In particular, we construct
plug-in classifiers that can achieve not only the fast, but also the
{\it super-fast} rates, i.e., the rates faster than $n^{-1}$. We
establish minimax lower bounds showing that the obtained rates
cannot be improved.
\end{abstract}



\noindent {\bf Short title.} Fast Rates for Plug-in Classifiers


\section{Introduction}

Let $(X,Y)$ be a random couple taking values in $\calZ \triangleq
\dsR^d\times\{0,1\}$ with joint distribution $P$. We regard $X\in
\dsR^d$ as a vector of features corresponding to an object and $Y\in
\{0,1\}$ as a label indicating that the object belongs to one of the
two classes. Consider the sample $(X_1, Y_1) , \ldots , (X_n , Y_n)
$, where $(X_i , Y_i)$ are independent copies of $(X,Y)$. We denote
by $\Pn$ the product probability measure according to which the
sample is distributed, and by $\px$ the marginal distribution of
$X$.

The goal of a classification procedure is to predict the label $Y$
given the value of $X$, i.e., to provide a decision rule $f:\dsR^d
\to\{0,1\}$ which belongs to the set $\calFXY$ of all Borel
functions defined on $\dsR^d$ and taking values in $\{0,1\}$. The
performance of a decision rule $f$ is measured by the
misclassification error
    $$R(f) \triangleq P(Y \neq f(X) ).$$
The Bayes decision rule is a minimizer of the risk $R(f)$ over all
the decision rules $f\in \calFXY$, and one of such minimizers has
the form $f^*(X)=\ds1_{\{\eta(X) \geq \demi\}}$ where
$\ds1_{\{\cdot\}}$ denotes the indicator function and $\eta(X)
\triangleq P(Y=1|X)$ is the regression function of $Y$ on $X$ (here
$P(dY|X)$ is a regular conditional probability which we will use in
the following without further mention).

An empirical decision rule (a classifier) is a random mapping $\hf :
\calZ^n \rightarrow \calFXY$ measurable w.r.t. the sample. Its
accuracy can be characterized by the excess risk
$$
{\cal E}(\hf) = \dsE R( \hf ) - R( f^* )
$$
where $\dsE$ is the sign of expectation. A key problem in
classification is to construct classifiers with small excess risk
for sufficiently large $n$ [cf. Devroye, Gy\"orfi and Lugosi (1996),
Vapnik (1998)]. Optimal classifiers can be defined as those having
the best possible rate of convergence of ${\cal E}(\hf)$ to 0, as
$n\to\infty$. Of course, this rate, and thus the optimal classifier,
depend on the assumptions on the joint distribution of $(X,Y)$. A
standard way to define optimal classifiers is to introduce a class
of joint distributions of $(X,Y)$ and to declare $\hf$ optimal if it
achieves the best rate of convergence in a minimax sense on this
class.

Two types of assumptions on the joint distribution of $(X,Y)$ are
commonly used: complexity assumptions and margin assumptions.

{\it Complexity assumptions} are stated in two possible ways. First
of them is to suppose that the regression function $\eta$ is smooth
enough or, more generally, belongs to a class of functions $\Sigma$
having a suitably bounded $\varepsilon$-entropy. This is called a
{\it complexity assumption on the regression function} (CAR). Most
commonly it is of the following form.

{\bf Assumption (CAR).} {\it The regression function $\eta$ belongs
to class $\Sigma$ of functions on $\dsR^d$ such that
$$
{\cal H}(\varepsilon,\Sigma,L_p) \le A_*\varepsilon^{-\rho}, \quad
\forall \varepsilon>0,
$$
with some constants $\rho>0$, $A_*>0$. Here ${\cal
H}(\varepsilon,\Sigma,L_p)$ denotes the $\varepsilon$-entropy of the
set $\Sigma$ w.r.t. an $L_p$ norm with some $1\le p\le \infty$.}

At this stage of discussion we do not identify precisely the value
of $p$ for the $L_p$ norm in Assumption (CAR), nor the measure with
respect to which this norm is defined. Examples will be given later.
If $\Sigma$ is a class of smooth functions with smoothness parameter
$\beta$ on a compact in $\dsR^d$, for example, a H\"older class, as
described below, a typical value of $\rho$ in Assumption (CAR) is
$\rho= d/\beta$.

Assumption (CAR) is well adapted for the study of {\it plug-in
rules}, i.e. of the classifiers having the form
\begin{equation}\label{pi}
\hf^{PI}(X) = \ds1_{\{\heta(X) \geq \demi\}}
\end{equation}
where $\heta$ is a nonparametric estimator of the function $\eta$.
Indeed, Assumption (CAR) typically reads as a smoothness assumption
on $\eta$ implying that a good nonparametric estimator (kernel,
local polynomial, orthogonal series or other) $\heta$ converges with
some rate to the regression function $\eta$, as $n\to\infty$. In
turn, closeness of $\heta$ to $\eta$ implies closeness of $\hf$ to
$f$: for any plug-in classifier $\hf^{PI}$ we have
\begin{equation}\label{pi1}
\dsE R( \hf^{PI} ) - R( f^* )\le 2 \dsE\int|\heta(x)-\eta(x)|\px
(dx)
\end{equation}
(cf. Devroye, Gy\"orfi and Lugosi (1996), Theorem 2.2). For various
types of estimators $\heta$ and under rather general assumptions it
can be shown that, if (CAR) holds, the RHS of (\ref{pi1}) is
uniformly of the order $n^{-1/(2+\rho)}$, and thus
\begin{equation}\label{pi2}
\sup_{P:\eta\in\Sigma} {\cal E}( \hf^{PI} ) = O(n^{-1/(2+\rho)}),
\quad n\to\infty,
\end{equation}
[cf. Yang (1999)]. In particular, if $\rho= d/\beta$ (which
corresponds to a class of smooth functions with smoothness parameter
$\beta$), we get
\begin{equation}\label{pi3}
\sup_{P:\eta\in\Sigma} {\cal E}( \hf^{PI} ) =
O(n^{-\beta/(2\beta+d)}), \quad n\to\infty.
\end{equation}
Note that (\ref{pi3}) can be easily deduced from (\ref{pi1}) and
standard results on the $L_1$ or $L_2$ convergence rates of usual
nonparametric regression estimators on $\beta$-smoothness classes
$\Sigma$. The rates in (\ref{pi2}), (\ref{pi3}) are quite slow,
always slower than $n^{-1/2}$. In (\ref{pi3}) they deteriorate
dramatically as the dimension $d$ increases. Moreover, Yang (1999)
shows that, under general assumptions, the bound (\ref{pi3}) cannot
be improved in a minimax sense. These results raised some pessimism
about the plug-in rules.

The second way to describe complexity is to introduce a structure on
the class of possible decision sets $G^* = \{x:f^*(x)=1\} =
\{x:\eta(x)\ge 1/2\}$ rather than on that of regression functions
$\eta$. A standard {\it complexity assumption on the decision set}
(CAD) is the following.

{\bf Assumption (CAD).} {\it The decision set $G^*$ belongs to a
class $\cal G$ of subsets of $\dsR^d$ such that
$$
{\cal H}(\varepsilon,{\cal G}, d_\triangle) \le
A_*\varepsilon^{-\rho}, \quad \forall \varepsilon>0,
$$
with some constants $\rho>0$, $A_*>0$. Here ${\cal
H}(\varepsilon,{\cal G},d_\triangle)$ denotes the
$\varepsilon$-entropy of the class $\cal G$ w.r.t. the measure of
symmetric difference pseudo-distance between sets defined by
$d_\triangle(G,G')= \px (G\triangle G')$ for two measurable subsets
$G$ and $G'$ in $\dsR^d$.}

The parameter $\rho$ in Assumption (CAD) typically characterizes the
smoothness of the boundary of $G^*$ [cf. Tsybakov (2004a)]. Note
that, in general, there is no connection between Assumptions (CAR)
and (CAD). Indeed, the fact that $G^*$ has a smooth boundary does
not imply that $\eta$ is smooth, and vice versa. The values of
$\rho$ closer to 0 correspond to smoother boundaries (less complex
sets $G^*$). As a limit case when $\rho\to 0$ one can consider the
Vapnik-Chervonenkis classes (VC-classes) for which the
$\varepsilon$-entropy is logarithmic in $1/\varepsilon$.

Assumption (CAD) is suited for the study of empirical risk
minimization (ERM) type classifiers introduced by Vapnik and
Chervonenkis (1974), see also Devroye, Gy\"orfi and Lugosi (1996),
Vapnik (1998). As shown in Tsybakov (2004a), for every $0<\rho<1$
there exist ERM classifiers $\hf^{ERM}$ such that, under Assumption
(CAD),
\begin{equation}\label{pi4}
\sup_{P:G^*\in{\cal G}} {\cal E}( \hf^{ERM} ) = O(n^{-1/2}), \quad
n\to\infty.
\end{equation}
The rate of convergence in (\ref{pi4}) is better than that for
plug-in rules, cf. (\ref{pi2}) -- (\ref{pi3}), and it does not
depend on $\rho$ (respectively, on the dimension $d$). Note that the
comparison between (\ref{pi4}) and (\ref{pi2}) -- (\ref{pi3}) is not
quite legitimate, because there is no inclusion between classes of
joint distributions $P$ of $(X,Y)$ satisfying Assumptions (CAR) and
(CAD). Nevertheless, such a comparison have been often interpreted
as an argument in disfavor of the plug-in rules. Indeed, Yang's
lower bound shows that the $n^{-1/2}$ rate cannot be attained under
Assumption (CAD) suited for the plug-in rules. Recently, advantages
of the ERM type classifiers, including penalized ERM methods, have
been further confirmed by the fact that, under the margin (or low
noise) assumption, they can attain {\it fast rates of convergence},
i.e. the rates that are faster than $n^{-1/2}$ [Mammen and Tsybakov
(1999), Tsybakov (2004a), Massart and N\'ed\'elec (2003), Tsybakov
and van de Geer (2005), Koltchinskii (2005), Audibert (2004)].

The {\it margin assumption} (or low noise assumption) is stated as
follows.

{\bf Assumption (MA).} {\it There exist constants $C_0>0$ and
$\alpha\ge 0$ such that
\begin{equation}\label{ma}
    \px\big(0<| \eta( X ) - 1/2 | \leq t \big) \leq C_0 t^{\alpha},
    \qquad \forall \  t>0.
\end{equation}
}

The case $\alpha=0$ is trivial (no assumption) and is included for
notational convenience. Assumption (MA) provides a useful
characterization of the behavior of regression function $\eta$ in a
vicinity of the level $\eta=1/2$ which turns out to be crucial for
convergence of classifiers (for more discussion of the margin
assumption see Tsybakov (2004a)). The main point is that, under
(MA), fast classification rates up to $n^{-1}$ are achievable. In
particular, for every $0<\rho<1$ and $\alpha>0$ there exist ERM type
classifiers $\hf^{ERM}$ such that
\begin{equation}\label{pi5}
\sup_{P: (CAD),(MA)} {\cal E}( \hf^{ERM} ) =
O(n^{-\frac{1+\alpha}{2+\alpha +\alpha\rho}}), \quad n\to\infty,
\end{equation}
where $\sup_{P: (CAD),(MA)}$ denotes the supremum over all joint
distributions $P$ of $(X,Y)$ satisfying Assumptions (CAD) and (MA).
The RHS of (\ref{pi5}) can be arbitrarily close to $O(n^{-1})$ for
large $\alpha$ and small $\rho$. Result (\ref{pi5}) for direct ERM
classifiers on $\varepsilon$-nets is proved by Tsybakov (2004a), and
for some other ERM type classifiers by Tsybakov and van de Geer
(2005), Koltchinskii (2005) and Audibert (2004) (in some of these
papers the rate of convergence (\ref{pi5}) is obtained with an extra
log-factor).

Comparison of (\ref{pi4}) and (\ref{pi5}) with (\ref{pi1}) and
(\ref{pi2}) seems to confirm the conjecture that the plug-in
classifiers are inferior to the ERM type ones. The main message of
the present paper is to disprove this conjecture. We will show that
there exist plug-in rules that converge with fast rates, and even
with {\it super-fast rates}, i.e. faster than $n^{-1}$ under the
margin assumption (MA). The basic idea of the proof  is to use
exponential inequalities for the regression estimator $\heta$ (see
Section 3 below) or the convergence results in the $L_\infty$ norm
(see Section 5), rather than the usual $L_1$ or $L_2$ norm
convergence of $\heta$, as previously described (cf. (\ref{pi1})).
We do not know whether the super-fast rates are attainable for ERM
rules or, more precisely, under Assumption (CAD) which serves for
the study of the ERM type rules. It is important to note that our
results on fast rates cover more general setting than just
classification with plug-in rules. These are rather results about
{\it classification in the regression complexity context under the
margin assumption}. In particular, we establish minimax lower bounds
valid for all classifiers, and we construct a ``hybrid" plug-in/ ERM
procedure (ERM based on a grid on a set regression functions $\eta$)
that achieves optimality. Thus, the point is mainly not about the
type of procedure (plug-in or ERM) but about the type of complexity
assumption (on the regression function (CAR) or on the decision set
(CAD)) that should be natural to impose. Assumption (CAR) on the
regression function arises in a natural way in the analysis of
several practical procedures of plug-in type, such as boosting and
SVM [cf. Blanchard, Lugosi and Vayatis (2003), Bartlett, Jordan and
McAuliffe (2003), Scovel and Steinwart (2003), Blanchard, Bousquet
and Massart (2004), Tarigan and van de Geer (2004)]. These
procedures are now intensively studied but, to our knowledge, only
suboptimal rates of convergence have been proved in the regression
complexity context under the margin assumption. The results in
Section 4 point out this fact (see also Section 5), and establish
the best achievable rates of classification that those procedures
should expectedly attain.

\section{Notation and definitions}

In this section we introduce some notation, definitions and basic
facts that will be used in the paper.

We denote by $C,C_1,C_2,\dots $ positive constants whose values may
differ from line to line. The symbols $\dsP$ and $\dsE$ stand for
generic probability and expectation signs, and $E_X$ is the
expectation w.r.t. the marginal distribution $\px$. We denote by
$\calB(x,r)$ the closed Euclidean ball in $\dsR^d$ centered at
$x\in\dsR^d$ and of radius $r>0$.

For any multi-index $s=(s_1,\dots,s_d)\in\dsN^d$ and any
$x=(x_1,\dots,x_d)\in\dsR^d$, we define $|s| = \sum_{i=1}^d s_i$,
$s!=s_1!\dots s_d!$, $x^s=x_1^{s_1} \dots x_d^{s_d}$ and $\|x\|
\triangleq (x_1^2+\cdots+x_d^2)^{1/2}$. Let $D^s$ denote the
differential operator
    $D^s\triangleq \frac{\partial^{s_1+ \cdots +s_d}}{\partial x_1^{s_1}
        \cdots \partial x_d^{s_d}}.$

Let $\beta >0$. Denote by $\lbeta$ the maximal integer that is
strictly less than $\beta$. For any $x\in\dsR^d$ and any $\lbeta$
times continuously differentiable real valued function $g$ on
$\dsR^d$, we denote by $g_x$ its Taylor polynomial of degree
$\lbeta$ at point $x$:
    $$g_x(x') \triangleq \sum_{|s|\le \lbeta}
    \frac{(x'-x)^s}{s!} D^s g(x).$$

Let $L>0$. The $\big(\beta,L,\dsR^d\big)$-{\it H\"{o}lder class} of
functions, denoted $\Sigma( \beta , L , \dsR^d )$, is defined as the
set of functions $g:\dsR^d\to \dsR$ that are $\lbeta$ times
continuously differentiable and satisfy, for any $x,x'\in\dsR^d$ ,
the inequality
        $$|g(x') - g_x(x')| \le L \|x-x'\|^{\beta}.$$

Fix some constants $c_0,r_0>0$. We will say that a Lebesgue
measurable set $A\subset \dsR^d$ is {\it $(c_0,r_0)$-regular} if
\begin{equation}\label{reg}
        \lam\big[A\cap\calB(x,r)\big] \ge c_0 \lam\big[\calB(x,r)\big],
        \quad \forall \ 0 < r \le r_0, \ \ \forall \ x \in A,
 \end{equation}
where $\lam[S]$ stands for the Lebesgue measure of $S\subset\dsR^d$.
 To illustrate this definition, consider the following example. Let
 $d\ge 2$. Then the set $A=\big\{ x=(x_1,\dots,x_d)\in\dsR^d: \sum_{j=1}^d
|x_j|^q \le 1 \big\}$ is $(c_0,r_0)$-regular with some $c_0,r_0>0$
for $q\ge 1$, and there are no $c_0,r_0>0$ such that $A$ is
$(c_0,r_0)$-regular for $0<q<1$.

Introduce now two assumptions on the marginal distribution $\px$
that will be used in the sequel.

\begin{defn} Fix $0< c_0,r_0,\mu_{\max}<\infty$ and a compact
$\calC\subset \dsR^d$. We say that the {\bf mild density assumption}
is satisfied if
    the
    marginal distribution $\px$ is supported on a compact $(c_0,r_0)$-regular set
    $A\subseteq\calC$
    and
has a uniformly bounded density $\mu$ w.r.t. the Lebesgue
    measure:
    \ $\mu(x)\le \mu_{\max}, \forall \ x\in A$.
\end{defn}

\begin{defn} Fix some constants $c_0,r_0>0$ and
$0<\mu_{\min} <\mu_{\max}<\infty$ and a compact $\calC\subset
\dsR^d$. We say that the {\bf strong density assumption} is
satisfied if
    the
    marginal distribution $\px$ is supported on a compact $(c_0,r_0)$-regular set
    $A\subseteq\calC$
    and has a density $\mu$ w.r.t. the Lebesgue measure bounded away from zero and infinity
    on $A$:
$$
            \mu_{\min} \le \mu(x) \le \mu_{\max}
            \quad \text{for} \ \ x\in A, \ \ \text{and} \ \
            \mu(x) = 0 \ \ \text{otherwise}.
$$
\end{defn}

We finally recall some notions related to locally polynomial
estimators.

\begin{defn}
For $h>0$, $x\in \dsR^d$, for an integer $l\ge0$ and a function
$K:\dsR^d\to \dsR_+$, denote by $\hth_x$ a polynomial on $\dsR^d$ of
degree $l$ which minimizes
\begin{equation}\label{lpe}
    \sum_{i=1}^n \big[Y_i - \hth_x(X_i-x)\big]^2 \kxi.
\end{equation}
 The {\bf locally polynomial estimator $\heta^{LP}(x)$ of order $l$},
or \lpl estimator, of the value $\eta(x)$ of the regression function
at point $x$ is defined by:
    $\heta^{LP}(x) \triangleq \hth_x(0)$ if $\hth_x$ is the
unique minimizer of (\ref{lpe}) and $\heta^{LP}(x) \triangleq 0$
otherwise. The value $h$ is called the bandwidth and the function
$K$ is called the kernel of the \lpl estimator.
\end{defn}

Let $T_s$ denote the coefficients of $\hth_x$ indexed by multi-index
$s\in\dsN^d$:
    $\hth_x(u) = \sum_{|s|\le l} T_s u^s.$
Introduce the vectors
    $T \triangleq \big(T_s\big)_{|s|\le l}$,
    $V \triangleq \big(V_s\big)_{|s|\le l}$ where
    \begarlab{eq:v}
    V_s \triangleq \sum_{i=1}^n Y_i (X_i-x)^s \kxi,
    \endarlab
    $U(u) \triangleq \big(u^s\big)_{|s|\le l}$
and the matrix $Q \triangleq \big(Q_{s_1,s_2}\big)_{|s_1|,|s_2|\le
l}$ where
    \begarlab{eq:q}
    Q_{s_1,s_2} \triangleq \sum_{i=1}^n (X_i-x)^{s_1+s_2} \kxi.
    \endarlab
The following result is straightforward (cf. Section 1.7 in Tsybakov
(2004b) where the case $d=1$ is considered).
\begin{prop}
If the matrix $Q$ is positive definite, there exists a unique
polynomial on $\dsR^d$ of degree $l$ minimizing (\ref{lpe}). Its
vector of coefficients is given by
    $T = Q^{-1} V$
and the corresponding \lpl regression function estimator has the
form
    $$\heta^{LP}(x) =U^T(0)Q^{-1}V= \sum_{i=1}^n Y_i \kxi U^T(0)Q^{-1}U(X_i-x).$$
\end{prop}

\section{Fast rates for plug-in rules: the
strong density assumption}

We first state a general result showing how the rates of convergence
of plug-in classifiers can be deduced from exponential inequalities
for the corresponding regression estimators.

In the sequel, for an estimator $\heta$ of $\eta$, we write
$$\dsP\big( \big|\heta(X) - \eta(X)| \ge
\delta \big)\triangleq  \int \Pn
        \big( \big|\heta(x) - \eta(x)| \ge \delta \big) \px(dx), \quad
        \forall \ \delta>0,
        $$
i.e., we consider the probability taken with respect to the
distribution of the sample $(X_1,Y_1,\dots X_n,Y_n)$ and the
distribution of the input $X$.

\begin{thm} \label{th:plugin}
Let $\heta$ be an estimator of the regression function $\eta$ and
$\calP$ a set of probability distributions on $\calZ$ such that for
some constants $C_1>0,$ $C_2>0$, for some positive sequence $a_n$,
for $n\ge 1$ and any $\delta>0$, and for almost all $x$ w.r.t.
$\px$, we have
    \begin{equation} \label{eq:heta}
    \underset{P\in\calP}{\sup} \, \Pn\Big( \big|\heta(x) - \eta(x)| \ge \delta \Big)
        \le C_1 \exp\big( - C_2 a_n \delta^2 \big).
    \end{equation}
Consider the plug-in classifier $\hf= \ds1_{\{\heta \geq
\demi\}}$. If all the distributions $P\in\calP$ satisfy the margin
assumption (MA), we have
    $$\underset{P\in\calP}{\sup} \Big\{ \dsE
        R( \hf ) - R( f^* ) \Big\} \le C a_n^{-\frac{1+\alpha}{2}}$$
for $n\ge 1$ with some constant $C>0$ depending only on $\alpha$,
$C_0$, $C_1$ and $C_2$.
\end{thm}

\begin{proof}
Consider the sets $A_j\subset\dsR^d, j=1,2,\dots,$ defined as
    \begar
    A_0 & \triangleq & \big\{ x\in\dsR^d:
    0<|\eta(x)-\demi| \le \delta \big\},\\
    A_j & \triangleq & \big\{ x\in\dsR^d: 2^{j-1} \delta < |\eta(x)-\demi| \le 2^j \delta \big\}, \qquad \text{ for } j\ge 1.
    \endar
For any $\delta>0$, we may write
    \begin{equation} \begin{array}{lll} \label{eq:peel}
    \dsE R(\hf) - R( f^* ) & = & \dsE\big(|2\eta(X)-1|
        \ds1_{\{\hf(X) \neq f^*(X)\}}\big)\\
    & = & \sum_{j=0}^\infty \dsE\big(|2\eta(X)-1| \ds1_{\{\hf(X) \neq f^*(X)\}}\ds1_{\{X\in A_j\}}\big) \\
    & \leq & 2\delta \px\big( 0<|\eta(X)-\demi|\le \delta\big)\\
    & & \qquad
        + \sum_{j\ge 1} \dsE\big(|2\eta(X)-1|
        \ds1_{\{\hf(X) \neq f^*(X)\}}\ds1_{\{X\in A_j\}} \big).
    \end{array} \end{equation}
On the event $\{\hf \neq f^*\}$ we have $|\eta - \demi|\le
|\heta-\eta|$. So, for any $j\ge 1$, we get
    \begar
    \dsE\big(|2\eta(X)-1|
        \ds1_{\{\hf(X) \neq f^*(X)\}}\ds1_{\{X\in A_j\}} \big)\\
    \qquad\qquad\qquad
        \le 2^{j+1} \delta \, \dsE\big[ \ds1_{\{|\heta(X)-\eta(X)| \ge 2^{j-1} \delta\}}
        \ds1_{\{0<|\eta(X)-\demi|\le 2^j \delta\}} \big]\\
    \qquad\qquad\qquad
        \le 2^{j+1} \delta \, E_X \Big[
        \Pn\big(|\heta(X)-\eta(X)| \ge 2^{j-1} \delta \big)
        \ds1_{\{0<|\eta(X)-\demi|\le 2^j \delta\}}\Big]\\
    \qquad\qquad\qquad
        \le C_1 2^{j+1} \delta \exp\Big( -C_2 a_n (2^{j-1}\delta)^2\Big)
        \px\Big( 0<|\eta(X)-\demi|\le 2^j \delta \Big)\\
    \qquad\qquad\qquad
        \le 2 C_1 C_0 2^{j(1+\alpha)} \delta^{1+\alpha}
        \exp\big( -C_2 a_n (2^{j-1}\delta)^2\big)\\
    \endar
where in the last inequality we used Assumption (MA). Now, from
inequality \eqref{eq:peel}, taking $\delta=a_n^{-1/2}$ and using
Assumption (MA) to bound the first term of the right hand side of
\eqref{eq:peel}, we get
    \begar
    \dsE R(\hf) - R( f^* ) & \le & 2 C_0 a_n^{-\frac{1+\alpha}{2}}
        + C a_n^{-\frac{1+\alpha}{2}} \sum_{j\ge 2} 2^{j(1+\alpha)} \exp\big( -C_2 2^{2j-2}\big)\\
    & \le & C a_n^{-\frac{1+\alpha}{2}}.
    \endar
\end{proof}

Inequality \eqref{eq:heta} is crucial to obtain the above result.
This inequality holds true for various types of estimators and
various sets of probability distributions $\calP$. Here we focus on
a standard case where $\eta$ belongs to the H\"older class
$\Sigma(\beta,L,\dsR^d)$ and the marginal law of $X$ satisfies the
strong density assumption. We are going to show that in this case
there exist estimators satisfying inequality \eqref{eq:heta} with
$a_n = n^{\frac{2 \beta}{2 \beta + d}}$. These can be, for example,
locally polynomial estimators. Specifically, assume from now on that
$K$ is a kernel satisfying
\begin{eqnarray}
&&\exists c>0:  \quad K(x) \ge c \ds1_{\{\|x\|\le c\}}, \quad
\forall x\in\dsR^d, \label{K1} \\
&&\int_{\dsR^d} K(u) du = 1,\label{K2}\\
&&\int_{\dsR^d} \big(1+\|u\|^{4\beta}\big) K^2(u) du < \infty,\label{K3}\\
&&\sup_{u\in\dsR^d} \big(1+\|u\|^{2\beta}\big) K(u) <
\infty.\label{K4}
\end{eqnarray}
Let $h>0$, and consider the matrix $\barO \triangleq
\big(\barO_{s_1,s_2}\big)_{|s_1|,|s_2|\le \lbeta}$ where
    $\barO_{s_1,s_2} = \frac{1}{nh^d} \sum_{i=1}^n \xih^{s_1+s_2} \kxi.$
Define the regression function estimator $\heta^*$  as follows. If
the smallest eigenvalue of the matrix $\barO$ is greater than $(\log
n)^{-1}$ we set $\heta^*(x)$ equal to the projection of
$\heta^{LP}(x)$ on the interval $[0,1]$, where $\heta^{LP}(x)$ is
the LP($\lbeta$) estimator with a bandwidth $h>0$ and a kernel $K$
satisfying (\ref{K1}) -- (\ref{K4}). If the smallest eigenvalue of
$\barO$ is less than $(\log n)^{-1}$ we set $\heta^*(x) = 0$.

\begin{thm} \label{th:kernel} Let $\calP$ be a class of probability
distributions $P$ on $\calZ$ such that the regression function
$\eta$ belongs to the H\"older class $\Sigma(\beta,L,\dsR^d)$ and
the marginal law of $X$ satisfies the strong density assumption.
Then there exist constants $C_1, C_2, C_3>0$ such that for any $0 <
h \le r_0/c$, any $C_3 h^\beta<\delta$ and any $n\ge1$ the estimator
$\heta^*$ satisfies
\begin{equation}\label{expon1}
    \underset{P\in\calP}{\sup} \,
        \Pn\Big( \big|\heta^*(x) - \eta(x)\big| \ge \delta \Big)
        \le C_1 \exp\big( - C_2 n h^d \delta^2 \big)
         \end{equation}
for almost all $x$ w.r.t. $\px$. As a consequence, there exist $C_1,
C_2>0$ such that for $h=n^{-\frac{1}{2\beta+d}}$ and any $\delta>0$,
$n\ge1$ we have
\begin{equation}\label{expon2}
\underset{P\in\calP}{\sup} \,
        \Pn\Big( \big|\heta^*(x) - \eta(x)\big| \ge \delta \Big)
        \le C_1 \exp\big( - C_2 n^{\frac{2\beta}{2\beta+d}}
        \delta^2 \big)
        \end{equation}
for almost all $x$ w.r.t. $\px$. The constants $C_1, C_2, C_3$
depend only on $\beta$, $d$, $L$, $c_0$, $r_0$, $\mu_{\min},
\mu_{\max}$, and on the kernel $K$.
\end{thm}

\begin{proof}
See Section \ref{sec:proofkernel}.
\end{proof}

\begin{rem}
We have chosen here the LP estimators of $\eta$ because for them the
exponential inequality \eqref{eq:heta} holds without additional
smoothness conditions on the marginal density of $X$. For other
popular regression estimators, such as kernel or orthogonal series
ones, similar inequality can be also proved if we assume that the
marginal density of $X$ is as smooth as the regression function.
\end{rem}

 \begin{defn}\label{calPS}
 For a fixed parameter $\alpha\ge 0$, fixed positive parameters $c_0,r_0,C_0,\beta, L,
\mu_{\max} > \mu_{\min} > 0$ and a fixed compact ${\cal C}\subset \dsR^d$, let $\calPS$ denote the class of all
probability distributions $P$ on $\cal Z$ such that
    \begin{enumerate}
    \item[(i)] the margin assumption (MA) is satisfied,
    \item[(ii)] the regression function $\eta$ belongs to the
    H\"{o}lder class $\Sigma(\beta,L,\dsR^d)$,
    \item[(iii)] the strong density assumption on $\px$ is satisfied.
    \end{enumerate}
\end{defn}

Theorem \ref{th:plugin} and (\ref{expon2}) immediately imply the
next result.

\begin{thm} \label{th:plugin1}
For any $n\ge1$ the excess risk of the plug-in classifier $\hf^* =
\ds1_{\{\heta^* \geq \demi\}}$ with bandwidth
$h=n^{-\frac{1}{2\beta+d}}$ satisfies
    $$\underset{P\in\calPS}{\sup} \Big\{ \dsE
        R( \hf^* ) - R( f^* ) \Big\}
        \le C n^{-\frac{\beta(1+\alpha)}{2\beta+d}}$$
where the constant $C>0$ depends only on $\alpha$, $C_0$, $C_1$ and
$C_2$.
\end{thm}

For $\alpha \beta
> d/2$ the convergence rate $n^{- \frac{\beta(1+\alpha)}{2\beta+d}}$
obtained in Theorem \ref{th:plugin1} is a {\it fast rate}, i.e., it
is faster than $n^{-1/2}$. Furthermore, it is a {\it super-fast
rate} (i.e., is faster than $n^{-1}$) for $(\alpha-1)\beta > d$. We must
note that if this condition is satisfied, the class $\calPS$ is
rather poor, and thus super-fast rates can occur only for very
particular joint distributions of $(X,Y)$. Intuitively, this is
clear. Indeed, to have a very smooth regression function $\eta$
(i.e., very large $\beta$) implies that when $\eta$ hits the level
$1/2$, it cannot ``take off'' from this level too abruptly. As a
consequence, when the density of the distribution $\px$ is bounded
away from $0$ at a vicinity of the hitting point, the margin
assumption cannot be satisfied for large $\alpha$ since this
assumption puts an upper bound on the ``time spent'' by the
regression function near $1/2$. So, $\alpha$ and $\beta$ cannot be
simultaneously very large. It can be shown that the cases of ``too
large" and ``not too large" $(\alpha,\beta)$ are essentially
described by the condition $(\alpha-1)\beta > d$.

To be more precise, observe first that $\calPS$ is not empty for
$(\alpha-1)\beta > d$, so that the super-fast rates can effectively
occur. Examples of laws $P\in\calPS$ under this condition can be
easily given, such as the one with $\px$ equal to the uniform
distribution on a ball centered at 0 in $\dsR^d$, and the regression
function defined by $\eta(x)=1/2-C\|x\|^2$ with an appropriate
$C>0$. Clearly, $\eta$ belongs to H\"older classes with arbitrarily
large $\beta$ and Assumption (MA) is satisfied with $\alpha=d/2$.
Thus, for $d\ge 3$ and $\beta$ large enough super-fast rates can
occur. Note that in this example the decision set $\{x: \eta(x)\ge
1/2\}$ has the Lebesgue measure 0 in $\dsR^d$. It turns out that
this condition is necessary to achieve classification with
super-fast rates when the H\"older classes of regression functions
are considered.

To explain this and to have further insight into the problem of
super-fast rates, consider the following two questions:
\begin{itemize}
\item for which parameters $\alpha,$ $\beta$ and $d$ is there a
distribution $P\in\calPS$ such that the regression function
associated with $P$ hits\footnote{ A function $f: \dsR^d \rightarrow
\dsR$ is said to hit the level $a\in\dsR$ at $x_0\in\dsR^d$ if
and only if $f(x_0)=a$ and for any $r>0$ there exists
$x \in \calB(x_0,r)$ such that $f(x) \neq a$
.} $1/2$ in the support of $\px$?
\item  for which parameters $\alpha,$ $\beta$ and $d$ is there a
distribution $P\in\calPS$ such that the regression function
associated with $P$ crosses\footnote{ A function $f: \dsR^d
\rightarrow \dsR$ is said to cross the level $a\in\dsR$ at
$x_0\in\dsR^d$ if and only if for any $r>0$, there exists $x_-$ and $x_+$ in
$\calB(x_0,r)$ such that $f(x_-)<a$ and $f(x_+)>a$.} $1/2$ in the
interior of the support of $\px$?
\end{itemize}

The following result gives a precise description of the constraints
on $(\alpha, \beta)$ leading to possibility or impossibility of the
super-fast rates.

\begin{prop} \label{th:hitcross}
\begin{itemize}
\item If $\alpha(1 \wedge \beta) > d$, there is no distribution $P\in\calPS$
such that the regression function $\eta$ associated with $P$ hits $1/2$
in the interior of the support of $P_X$.
\item For any $\alpha,\beta>0$ and integer $d\ge\alpha(1 \wedge \beta)$, any
positive parameter $L$
and any compact $\calC\subset\dsR^d$ with non-empty interior,
for appropriate positive parameters $C_0,c_0,r_0,\mu_{\max}>\mu_{\min}>0$,
there are distributions $P\in\calPS$ such that the regression function $\eta$
associated with $P$ hits $1/2$ in the boundary of the support of $P_X$.
\item For any $\alpha,\beta>0$, any integer $d\ge 2\alpha$, any
positive parameter $L$ and any compact $\calC\subset\dsR^d$ with non-empty interior,
for appropriate positive parameters $C_0,c_0,r_0,\mu_{\max}>\mu_{\min}>0$,
there are distributions $P\in\calPS$ such that the regression function $\eta$ associated with $P$ hits $1/2$
in the interior of the support of $P_X$.
\item If $\alpha (1 \wedge \beta) > 1$ there is no distribution
$P\in\calPS$ such that the regression function $\eta$ associated
with $P$ crosses $1/2$ in the interior of the support of $\px$.
Conversely, for any $\alpha,\beta>0$ such that $\alpha (1 \wedge \beta) \le 1$,
any integer $d$, any positive parameter $L$
and any compact $\calC\subset\dsR^d$ with non-empty interior,
for appropriate positive parameters $C_0,c_0,r_0,\mu_{\max}>\mu_{\min}>0$,
there are distributions $P\in\calPS$ such that the regression function $\eta$ associated with $P$ crosses $1/2$
in the interior of the support of $P_X$.
\end{itemize}
\end{prop}

Note that the condition $\alpha (1 \wedge \beta) > 1$ appearing in
the last assertion is equivalent to
$\frac{\beta(1+\alpha)}{2\beta+d} >
\frac{(2\beta)\vee(\beta+1)}{2\beta+d}$, which is necessary to have
super-fast rates. As a consequence, in this context,
super-fast rates cannot occur when the regression function crosses
$1/2$ in the interior of the support. The third assertion of the proposition
shows that super-fast rates can occur with regression functions
hitting $1/2$ in the interior of the support of $\px$ provided that the regression
function is highly smooth and defined on a highly dimensional space and that
a strong margin assumption holds (i.e. $\alpha$ large).

\begin{proof}
See Section \ref{sec:proofhc}.
\end{proof}

The following lower bound shows optimality of the rate of
convergence for the H\"{o}lder classes obtained in Theorem
\ref{th:plugin1}.

\begin{thm} \label{th:lb}
Let $d\ge 1$ be an integer, and let $L,\beta,\alpha$ be positive
constants, such that $\alpha\beta\le d$. Then there exists a
constant $C>0$ such that for any $n\ge1$ and any classifier $\hf :
\calZ^n \rightarrow \calFXY$, we have
        $$\begin{array}{lll}
        \underset{P\in \calPS}{\sup} \big\{ \dsE R( \hf ) - R( f^* ) \big\}
            \geq C n^{- \frac{\beta(1+\alpha)}{2\beta+d}}.
        \end{array}$$
\end{thm}

\begin{proof}
See Section \ref{sec:prooflb}.
\end{proof}

Note that the lower bound of Theorem \ref{th:lb} does not cover the
case of super-fast rates ($(\alpha-1) \beta > d$).

Finally, we discuss the case where ``$\alpha=\infty$", which means
that there exists $t_0>0$ such that
\begin{equation}\label{infty}
\px \big(0<|\eta(X)-1/2|\le t_0\big) =0.
\end{equation}
This is a very favorable situation for classification. The rates of
convergence of the ERM type classifiers under (\ref{infty}) are, of
course, faster than under Assumption (MA) with $\alpha<\infty$ [cf.
Massart and N\'ed\'elec (2003)], but they are not faster than
$n^{-1}$. Indeed, Massart and N\'ed\'elec (2003) provide a lower
bound showing that, even if Assumption (CAD) is replaced by a very
strong assumption that the true decision set belongs to a VC-class
(note that both assumptions are naturally linked to the study the
ERM type classifiers), the best achievable rate is of the order
$(\log n)/n$. We show now that for the plug-in classifiers much
faster rates can be attained. Specifically, if the regression
function $\eta$ has some (arbitrarily low) H\"older smoothness
$\beta$ the rate of convergence can be exponential in $n$. To show
this, we first state a simple lemma which is valid for any plug-in
classifier $\hf$.
\begin{lem}\label{expc}
Let assumption (\ref{infty}) be satisfied, and let $\heta$ be an
estimator of the regression function $\eta$. Then for the plug-in
classifier $\hf = \ds1_{\{\heta \geq \demi\}}$ we have
$$\dsE
        R( \hf ) - R( f^* )
        \le \dsP
        \big(|\heta(X)-\eta(X)| > t_0 \big).$$
\end{lem}
\begin{proof}
Following the argument similar to the proof of Theorem
\ref{th:plugin} and using condition (\ref{infty}) we get
    \begin{equation*} \begin{array}{lll} \label{eq?}
    \dsE R(\hf) - R( f^* )
    & \leq & 2t_0 \px\big(0< |\eta(X)-1/2|\le t_0\big)\\
    & & \qquad
        + \ \dsE\big(|2\eta(X)-1|
        \ds1_{\{\hf(X) \neq f^*(X)\}}\ds1_{\{|\eta(X)-1/2|> t_0\}}
        \big)\\
     & = & \dsE\big(|2\eta(X)-1|
        \ds1_{\{\hf(X) \neq f^*(X)\}}\ds1_{\{|\eta(X)-1/2|> t_0\}}
        \big)\\
    & \leq & \dsP
        \big(|\heta(X)-\eta(X)| > t_0 \big).
    \end{array} \end{equation*}
\end{proof}

Lemma \ref{expc} and Theorem \ref{th:kernel} immediately imply that,
under assumption (\ref{infty}), the rate of convergence of the
plug-in classifier
 $\hf^* =
\ds1_{\{\heta^* \geq \demi\}}$ with a small enough fixed
(independent of $n$) bandwidth $h$ is exponential. To state the
result, we denote by $\calPSI$ the class of probability
distributions $P$ defined in the same way as $\calPS$, with the only
difference that in Definition \ref{calPS} the margin assumption (MA)
is replaced by condition (\ref{infty}).
\begin{prop} \label{prop:plugin2}
There exists a fixed (independent of $n$) $h>0$ such that for any
$n\ge1$ the excess risk of the plug-in classifier $\hf^* =
\ds1_{\{\heta^* \geq \demi\}}$ with bandwidth $h$ satisfies
    $$\underset{P\in\calPSI}{\sup} \Big\{ \dsE
        R( \hf^* ) - R( f^* ) \Big\}
        \le C_4 \exp(-C_5 n)$$
where the constants $C_4,C_5>0$ depend only on $t_0$, $\beta$, $d$,
$L$, $c_0$, $r_0$, $\mu_{\min}, \mu_{\max}$, and on the kernel $K$.
\end{prop}
\begin{proof}
Use Lemma \ref{expc}, choose $h>0$ such that $h<\min(r_0/c,
(t_0/C_3)^{1/\beta})$, and apply (\ref{expon1}) with $\delta=t_0$.
\end{proof}

Koltchinskii and Beznosova (2005) prove a result on exponential
rates for the plug-in classifier with some penalized regression
estimator in place of the locally polynomial one that we use here.
Their result is stated under a less general condition, in the sense
that they consider only the Lipschitz class of regression functions
$\eta$, while in Proposition \ref{prop:plugin2} the H\"older
smoothness $\beta$ can be arbitrarily close to 0. Note also that we
do not impose any complexity assumption on the decision set.
However, the class of distributions $\calPSI$ is quite restricted in
a different sense. Indeed, for such distributions condition
(\ref{infty}) should be compatible with the assumption that $\eta$
belongs to a H\"older class. A sufficient condition for that is the
existence of a band or a ``corridor" of zero $P_X$-measure
separating the sets $\{x:\eta(x)>1/2\}$ and $\{x:\eta(x)<1/2\}$. We
believe that this condition is close to the necessary one.

\section{Optimal learning rates without the strong density assumption}

In this section we show that if $\px$ does not admit a density
bounded away from zero on its support the rates of classification
are slower than those obtained in Section 3. In particular,
super-fast rates, i.e., the rates faster than $n^{-1}$, cannot be
achieved. Introduce the following class of probability
distributions.

\begin{defn}\label{def41}
For a fixed parameter $\alpha\ge 0$, fixed positive parameters $c_0,r_0,C_0,\beta,
L, \mu_{\max} > 0$ and a fixed compact ${\calC}\subset \dsR^d$, let $\calPS'$ denote the class of all
probability distributions $P$ on $\calZ$ such that
    \begin{enumerate}
    \item[(i)] the margin assumption (MA) is satisfied,
    \item[(ii)] the regression function $\eta$ belongs to the
    H\"{o}lder class $\Sigma(\beta,L,\dsR^d)$,
    \item[(iii)] the mild density assumption on $\px$ is satisfied.
    \end{enumerate}
\end{defn}

In this section we mainly assume that the distribution $P$ of
$(X,Y)$ belongs to $\calPS'$, but we also consider larger classes of
distributions satisfying the margin assumption (MA) and the
complexity assumption (CAR).

Clearly, $\calPS\subset \calPS'$. The only difference between
$\calPS'$ and $\calPS$ is that for $\calPS'$ the marginal density of
$X$ is not bounded away from zero. The optimal rates for $\calPS'$
are slower than for $\calPS$. Indeed, we have the following lower
bound for the excess risk.

\begin{thm} \label{th:lb2}
Let $d\ge 1$ be an integer, and let $L,\beta,\alpha$ be positive
constants. Then there exists a constant $C>0$ such that for any
$n\ge1$ and any classifier $\hf : \calZ^n \rightarrow \calFXY$ we
have
        $$\begin{array}{lll}
        \underset{P \in \calPS'}{\sup} \big\{ \dsE R( \hf ) - R( f^* ) \big\}
            \geq C n^{- \frac{(1+\alpha)\beta}{(2+\alpha)\beta+d}}.
        \end{array}$$
\end{thm}

\begin{proof}
See Section \ref{sec:prooflb}.
\end{proof}

In particular, when $\alpha = {d}/{\beta}$, we get slow convergence
rate $1/\sqrt{n}$, instead of the fast rate
$n^{-\frac{\beta+d}{2\beta+d}}$ obtained in Theorem \ref{th:plugin1}
under the strong density assumption. Nevertheless, the lower bound
can still approach $n^{-1}$, as the margin parameter $\alpha$ tends
to $\infty$.

We now show that the rate of convergence given in Theorem
\ref{th:lb2} is optimal in the sense that there exist estimators
that achieve this rate. This will be obtained as a consequence of a
general upper bound for the excess risk of classifiers over a larger
set $\calP$ of distributions than $\calPS'$.

Fix a Lebesgue measurable set $\calC\subset \dsR^d$ and a value
$1\le p\le \infty$. Let $\Sigma$ be a class of regression functions
$\eta$ on $\dsR^d$ such that Assumption (CAR) is satisfied where the
$\varepsilon$-entropy is taken w.r.t. the $L_p(\calC,\lam)$ norm
($\lam$ is the Lebesgue measure on $\dsR^d$). Then for every
$\varepsilon>0$ there exists an $\varepsilon$-net
$\calN_\varepsilon$ on $\Sigma$ w.r.t. this norm such that
$$
\log \big({\rm card} \, \calN_\ve \big)\le A'\ve^{-\rho},
$$
where $A'$ is a constant. Consider the empirical risk
$$
R_n(f) = \frac{1}{n} \sum_{i=1}^n \ds1_{\{f(X_i)\ne Y_i\}}, \quad \
f\in{\cal F},
$$
and set
$$
\ve_n = \ve_n(\alpha,\rho,p) \triangleq \left\{
\begin{array}{lcr}
n^{-\frac{1}{2+\alpha+\rho}}& \text{if} & p=\infty, \\
n^{-\frac{p+\alpha}{(2+\alpha)p+\rho(p+\alpha)}}& \text{if} & 1\le
p<\infty.
\end{array}
\right.
$$
Define a sieve estimator $\heta^S$ of the regression function $\eta$
by the relation
\begin{equation} \label{enet}
\heta^S \in {\rm
Argmin}_{\bar\eta\in{\calN}_{\ve_n}}R_n(f_{\bar\eta})
\end{equation}
where $f_{\bar\eta}(x)=\ds1_{\{\bar\eta(x) \geq 1/2\}}$, and
consider the classifier ${\hat f}_n^S=\ds1_{\{\heta^S \geq 1/2\}}$.
Note that ${\hat f}_n^S$ can be viewed as a ``hybrid" plug-in/ ERM
procedure: the ERM is performed on a set of plug-in rules
corresponding to a grid on the class of regression functions $\eta$.

\begin{thm} \label{th:ub3}
Let $\calP$ be a set of probability distributions $P$ on $\calZ$
such that
    \begin{enumerate}
    \item[(i)] the margin assumption (MA) is satisfied,
    \item[(ii)] the regression function $\eta$ belongs to
    a class $\Sigma$ which satisfies the complexity assumption (CAR)
    with the $\varepsilon$-entropy
taken w.r.t. the $L_p(\calC,\lam)$ norm for some $1\le p\le \infty$,
    \item[(iii)] for all $P\in\calP$ the supports of
    marginal distributions $\px$ are included in $\calC$.
    \end{enumerate}
Consider the classifier ${\hat f}_n^S=\ds1_{\{\heta^S \geq 1/2\}}$.
If $p=\infty$ for any $n\ge1$ we have
\begin{equation}
\label{t42.1}\underset{P\in\calP}{\sup} \Big\{ \dsE
        R( {\hat f}_n^S ) - R( f^* ) \Big\} \le
C n^{- \frac{1+\alpha}{2+\alpha+\rho}}.
\end{equation}
If $1\le p<\infty$ and, in addition, for all $P\in\calP$ the
marginal distributions $\px$ are absolutely continuous w.r.t. the
Lebesgue measure and their densities are uniformly bounded from
above by some constant $\mu_{\max}<\infty$, then for any $n\ge1$ we
have
\begin{equation} \label{t42.2}
\underset{P\in\calP}{\sup} \Big\{ \dsE
        R( {\hat f}_n^S ) - R( f^* ) \Big\} \le
 C n^{- \frac{(1+\alpha)p}{(2+\alpha)p+\rho(p+\alpha)}}.
\end{equation}
\end{thm}

\begin{proof}
See Section \ref{proof:ub2}.
\end{proof}

Theorem \ref{th:ub3} allows one to get fast classification rates
without any density assumption on $\px$. Namely, define the
following class of distributions $P$ of $(X,Y)$.

\begin{defn}\label{def42}
For fixed parameters $\alpha\ge 0$, $C_0>0,\beta>0, L> 0$, and for a
fixed compact $\calC\subset \dsR^d$, let $\calPS^{0}$ denote the
class of all probability distributions $P$ on $\calZ$ such that
    \begin{enumerate}
    \item[(i)] the margin assumption (MA) is satisfied,
    \item[(ii)] the regression function $\eta$ belongs to the
    H\"{o}lder class $\Sigma(\beta,L,\dsR^d)$,
    \item[(iii)] for all $P\in\calP$ the supports of
    marginal distributions $\px$ are included in $\calC$.
    \end{enumerate}
\end{defn}

If $\calC$ is a compact the estimates of $\varepsilon$-entropies of
H\"{o}lder classes $\Sigma(\beta,L,\dsR^d)$ in the
$L_\infty(\calC,\lam)$ norm can be obtained from Kolmogorov and
Tikhomorov~(1961), and they yield Assumption (CAR) with
$\rho=d/\beta$. Therefore, from (\ref{t42.1}) we easily get the
following upper bound.

\begin{thm} \label{th:ub2}
Let $d\ge 1$ be an integer, and let $L,\beta,\alpha$ be positive
constants. For any $n\ge 1$ the classifier ${\hat
f}_n^S=\ds1_{\{\heta^S \geq 1/2\}}$ defined by (\ref{enet}) with
$p=\infty$ satisfies
    $$\underset{P\in\calPS^{0}}{\sup} \Big\{ \dsE
        R( {\hat f}_n^S ) - R( f^* ) \Big\} \le C n^{- \frac{(1+\alpha)\beta}
        {(2+\alpha)\beta+d}}$$
with some constant $C>0$ depending only on $\alpha$, $\beta$, $d$,
$L$ and $C_0$.
\end{thm}

Since $\calPS'\subset \calPS^{0}$, Theorems \ref{th:lb} and
\ref{th:ub2} show that $n^{-
\frac{(1+\alpha)\beta}{(2+\alpha)\beta+d}}$ is optimal rate of
convergence of the excess risk on the class of distributions
$\calPS^{0}$.

\section{Comparison lemmas}

In this section we give some useful inequalities between the risks
of plug-in classifiers and the $L_p$ risks of the corresponding
regression estimators under the margin assumption (MA). These
inequalities will be helpful in the proofs. They also illustrate a
connection between the two complexity assumptions (CAR) and (CAD)
defined in the Introduction and allow one to compare our study of
plug-in estimators with that given by Yang (1999) who considered the
case $\alpha=0$ (no margin assumption),
as well as with the developments in Bartlett, Jordan and McAuliffe
(2003) and Blanchard, Lugosi and Vayatis (2003).

Throughout this section $\bar\eta$ is a Borel function on $\dsR^d$
and
$$
\bar f(x)=\ds1_{\{\bar\eta(x) \geq 1/2\}}.
$$
For $1\le p\le \infty$ we denote by $\|\cdot\|_p$ the $L_p(\dsR^d,
\px)$ norm. We first state some comparison inequalities for the
$L_\infty$ norm.

\begin{lem} \label{lem:comp}
For any distribution $P$ of $(X,Y)$ satisfying Assumption $(MA)$ we
have
\begin{equation}        \label{eq:comp1}
        R(\bar f  ) - R( f^* )
        \le 2C_0 \|\bar\eta - \eta\|_\infty^{1+\alpha},
        \end{equation}
and
\begin{equation}        \label{eq:comp2}
        \px \big(\bar f (X) \neq f^*(X) , \, \eta(X)\neq 1/2\big)
        \le C_0 \|\bar\eta - \eta\|_\infty^{\alpha}.
        \end{equation}
\end{lem}
\begin{proof} To show (\ref{eq:comp1}) note that
    \begar
    R(\bar f)-R(f^*) & =
    & \dsE\big( |2\eta(X)-1| \ds1_{\{\bar f(X)\neq f^*(X)\}} \big)\\
    & \le & 2 \dsE\big( |\eta(X)-\demi| \ds1_{0<\{|\eta(X)-\demi|
    \le |\eta(X)-\bar\eta(X)|\}} \big)\\
    & \le & 2 \|\eta-\bar\eta\|_\infty \px\big(0<|\eta(X)-\demi|
    \le \|\eta-\bar\eta\|_\infty\big)\\
    & \le & 2 C_0 \|\eta-\bar\eta\|_\infty ^{1+\alpha}.
    \endar
Similarly,
 \begar
\px\big(\bar f(X)\neq f^*(X),\, \eta(X)\neq 1/2\big)
    & \le & \px\big(0<|\eta(X)-\demi|\le |\eta(X)-\bar\eta(X)|\big)\\
    & \le & \px\big(0<|\eta(X)-\demi|\le \|\eta-\bar\eta\|_\infty\big)\\
    & \le & C_0 \|\eta-\bar\eta\|_\infty^\alpha.
 \endar
\end{proof}

\begin{rem}
Lemma \ref{lem:comp} offers an easy way to obtain the result of
Theorem \ref{th:plugin1} in a slightly less precise form, with an
extra logarithmic factor in the rate. In fact, under the strong
density assumption, there exist nonparametric estimators $\heta$
(for instance, suitably chosen locally polynomial estimators) such
that
$$
\dsE \Big(\|\heta-\eta\|_\infty^q \Big)\le C\left(\frac{\log
n}{n}\right)^{\frac{q\beta}{2\beta+d}}, \quad \forall \ q>0,
$$
uniformly over $\eta\in\Sigma(\beta,L,\dsR^d)$ [see, e.g., Stone
(1982)]. Taking here $q=1+\alpha$ and applying the comparison
inequality (\ref{eq:comp1}) we immediately get that the plug-in
classifier $ \hf=\ds1_{\{\heta \geq 1/2\}}$ has the excess risk
${\cal E}(\hf)$ of the order $\left(n/\log
n\right)^{-\beta(1+\alpha)/(2\beta+d)}$.
\end{rem}

Another immediate application of Lemma \ref{lem:comp} is to get
lower bounds on the risks of regression estimators in the $L_\infty$
norm from the corresponding lower bounds on the excess risks of
classifiers (cf. Theorems \ref{th:lb} and \ref{th:lb2}). But here
again we loose a logarithmic factor required for the best bounds.

We now consider the comparison inequalities for $L_p$ norms with
$1\le p <\infty$.
\begin{lem} \label{lem:sqrtg}
For any $1\le p <\infty$ and any distribution $P$ of $(X,Y)$
satisfying Assumption $(MA)$ with $\alpha>0$ we have
\begin{equation}        \label{eq:comp3}
        R(\bar f  ) - R( f^* )
        \le C_1(\alpha,p) \|\eta-\bar\eta\|_p^
        {\frac{p(1+\alpha)}{p+\alpha}},
        \end{equation}
and
\begin{equation}        \label{eq:comp4}
        \px \big(\bar f (X) \neq f^*(X) ,\, \eta(X)\neq 1/2\big)
        \le C_2(\alpha,p) \|\eta-\bar\eta\|_p^
        {\frac{p}{p+\alpha}},
        \end{equation}
where $C_1(\alpha,p)=2(\alpha+p)p^{-1} \big( \frac{p}{\alpha}
        \big)^{\frac{\alpha}{\alpha+p}} C_0^{\frac{p-1}{\alpha+p}}$,
$C_2(\alpha,p)=(\alpha+p)p^{-1} \big( \frac{p}{\alpha}
        \big)^{\frac{\alpha}{\alpha+p}} C_0^{\frac{p}{\alpha+p}}$.
        In particular,
\begin{equation}        \label{eq:sqrtg}
        R(\bar f  ) - R( f^* )
        \le C_1(\alpha,2) \left( \int [ \bar\eta(x) - \eta(x) ]^2 \px (dx)
            \right)^{\frac{1+\alpha}{2+\alpha}}.
        \end{equation}
\end{lem}

\begin{proof}
For any $t > 0$ we have
        \begin{eqnarray}
        \lefteqn{R( \bar f ) - R( f^* )} \notag\\
        & & = \dsE \big[ |2\eta(X)-1| \ds1_{\{\bar f(X) \neq f^*(X)\}}
        \big]\notag\\
        & & = 2 \dsE \big[ |\eta(X)-1/2|
        \ds1_{\{\bar f(X) \neq f^*(X)\}}  \ds1_{\{0<|\eta(X)-1/2| \leq t\}}
        \big] \notag\\
        & & \mbox{} + 2\dsE \big[ |\eta(X)-1/2|
        \ds1_{\{\bar f(X) \neq f^*(X)\}} \ds1_{\{|\eta(X)-1/2|>t\}} \big] \notag\\
        & & \leq 2 \dsE \big[ |\eta(X)-\bar\eta(X)|
        \ds1_{\{0<|\eta(X)-1/2| \leq t\}} \big]
            + 2 \dsE \big[ |\eta(X)-\bar\eta(X)| \ds1_{\{|\eta(X)-\bar\eta(X)|>t\}} \big] \notag\\
        & & \leq 2  \|\eta-\bar\eta\|_p
            \big[\px ( 0<|\eta(X)-1/2| \leq t )\big]^{\frac{p-1}{p}}
            + \frac{2\|\eta-\bar\eta\|_p^p}{ t^{p-1} }
            \label{eq:sqrtg2}
        \end{eqnarray}
by H\"older and Markov inequalities. So, for any $t > 0$,
introducing $E\triangleq\|\eta-\bar\eta\|_p$
and using Assumption (MA) to bound the probability in
(\ref{eq:sqrtg2}) we obtain
        $$R( \bar f ) - R( f^* ) \leq 2 \left(
        C_0^{\frac{p-1}{p}} \, t^{\frac{\alpha(p-1)}{p}}  E +
        \frac{E^p}{t^{p-1}}\right).$$
Minimizing in $t$ the RHS of this inequality we get
(\ref{eq:comp3}). Similarly,
\begin{eqnarray*} \px\big(\bar f(X)\neq
f^*(X),\, \eta(X)\neq 1/2\big)
    & \le & \px\big(0<|\eta(X)-1/2|\le t\big)+
            \px\big(|\eta(X)-\bar\eta(X)|>t)\\
    & \le & C_0 t^\alpha + \frac{\|\eta-\bar\eta\|_p^p}{t^p} \ ,
\end{eqnarray*}
and minimizing this bound in $t$ we obtain (\ref{eq:comp4}).
\end{proof}

If the regression function $\eta$ belongs to the H\"{o}lder class
$\Sigma\big(\beta,L,\dsR^d\big)$ there exist estimators $\heta$ such
that, uniformly over the class,
    \begin{equation}\label{sy}
    \dsE\Big\{ \big[\heta( X ) - \eta(X)\big]^2 \Big\}
    \le C n^{-\frac{2\beta}{2\beta+d}}\end{equation}
for some constant $C>0$. This has been shown by Stone (1982) under
the additional strong density assumption and by Yang (1999) with no
assumption on $\px$. Using (\ref{sy}) and (\ref{eq:sqrtg}) we get
that the excess risk of the corresponding plug-in classifier $
\hf=\ds1_{\{\heta \geq 1/2\}}$ admits a bound of the order
$n^{-\frac{2\beta}{2\beta+d}\frac{1+\alpha}{2+\alpha}}$ which is
suboptimal when $\alpha \neq 0$ (cf. Theorems \ref{th:ub3},
\ref{th:ub2}). In other words, under the margin assumption, Lemma
\ref{lem:sqrtg} is not the right tool to analyze the convergence
rate of plug-in classifiers. On the contrary, when no margin
assumption is imposed (i.e., $\alpha=0$ in our notation) inequality
(\ref{pi1}), which is a version of (\ref{eq:sqrtg}) for $\alpha=0$,
is precise enough to give the optimal rate of classification [Yang
(1999)].

Another way to obtain (\ref{eq:sqrtg}) is to use Bartlett, Jordan
and McAuliffe (2003): it is enough to apply their Theorem 10 with
(in their notation) $\phi(t)=(1-t)^2,\psi(t)=t^2$ and to note that
for this choice of $\phi$ we have $R_\phi (\bar\eta) - R_\phi^* =
\|\eta-\bar\eta\|_2^2$. Blanchard, Lugosi and Vayatis (2003) used
the result of Bartlett, Jordan and McAuliffe (2003) to prove fast
rates of the order $n^{-\frac{2(1+\alpha)}{3(2+\alpha)}}$ for a
boosting procedure over the class of regression functions $\eta$ of
bounded variation in dimension $d=1$. Note that the same rates can
be obtained for other plug-in classifiers using (\ref{eq:sqrtg}).
Indeed, if $\eta$ is of bounded variation, there exist estimators of
$\eta$ converging with the mean squared $L_2$ rate $n^{-2/3}$[cf.
van de Geer (2000), Gy\"{o}rfi et al. (2002)], and thus application
of (\ref{eq:sqrtg}) immediately yields the rate
$n^{-\frac{2(1+\alpha)}{3(2+\alpha)}}$ for the corresponding plug-in
rule. However, Theorem \ref{th:ub3} shows that this is not an
optimal rate (here again we observe that inequality (\ref{eq:sqrtg})
fails to establish the optimal properties of plug-in classifiers).
In fact, let $d=1$ and let the assumptions of Theorem \ref{th:ub3}
be satisfied, where instead of assumption (ii) we use its particular
instance: $\eta$ belongs to a class of functions on $[0,1]$ whose
total variation is bounded by a constant $L<\infty$. It follows from
Birman and Solomjak (1967) that Assumption (CAR) for this class is
satisfied with $\rho=1$ for any $1\le p<\infty$. Hence, we can apply
(\ref{t42.2}) of Theorem \ref{th:ub3} to find that
\begin{equation}\label{birman}
\underset{P\in\calP}{\sup} \Big\{ \dsE
        R( {\hat f}_n^S ) - R( f^* ) \Big\} \le
 C n^{- \frac{(1+\alpha)p}{(2+\alpha)p+(p+\alpha)}}
\end{equation}
for the corresponding class $\calP$. If $p>2$ (recall that the value
$p\in [1,\infty)$ is chosen by the statistician), the rate in
(\ref{birman}) is faster than $n^{-\frac{2(1+\alpha)}{3(2+\alpha)}}$
obtained under the same conditions by Blanchard, Lugosi and Vayatis
(2003).

\section{Proofs}

\subsection{Proof of Theorem \ref{th:kernel}} \label{sec:proofkernel}

Consider a distribution $P$ in $\calPS$. Let $A$ be the support of
$\px$. Fix $x\in A$ and $\delta>0$. Consider the matrix $B
\triangleq \big( B_{s_1,s_2}\big)_{|s_1|,|s_2|\le \lbeta}$ with
elements $B_{s_1,s_2} \triangleq \int_{\dsR^d} u^{s_1+s_2} K(u)
\mu(x+hu) du$. The smallest eigenvalue $\lam_{\barO}$ of $\barO$
satisfies
    \begarlab{eq:lamo}
    \lam_{\barO} & = & \min_{\|W\|=1} W^T \barO W\\
    & \ge & \min_{\|W\|=1} W^T B W + \min_{\|W\|=1} W^T (\barO - B) W\\
    & \ge & \min_{\|W\|=1} W^T B W - \sum_{|s_1|,|s_2|\le \lbeta} \big|
    \barO_{s_1,s_2}
        - B_{s_1,s_2}\big|.
    \endarlab
Let $A_n \triangleq \big\{ u \in \dsR^d : \|u\|\le c ; \ x+hu \in A
\big\}$
where $c$ is the constant appearing in (\ref{K1}). Using (\ref{K1}),
for any vector $W$ satisfying $\|W\|=1$, we obtain
    \begar
    W^T B W & = & \int_{\dsR^d} \big( \sum_{|s|\le \lbeta} W_s u^s \big)^2
        K(u) \mu(x+hu) du\\
    & \ge &
    c\mu_{\min} \int_{A_n} \big( \sum_{|s|\le \lbeta} W_s u^s \big)^2
    du.
    \endar
By assumption of the theorem, $ch \le r_0$. Since the support of the
marginal distribution is $(c_0,r_0)$-regular we get
    \begar
    \lam[A_n] \ge h^{-d} \lam\big[ \calB(x,ch) \cap A \big]
        \ge c_0 h^{-d} \lam\big[ \calB(x,ch) \big] \ge c_0 v_d c^d,
    \endar
where $v_d \triangleq \lam\big[\calB(0,1)\big]$ is the volume of the
unit ball and $c_0>0$ is the constant introduced in the definition
(\ref{reg}) of the $(c_0,r_0)$-regular set.

Let $\calA$ denote the class of all compact subsets of $\calB(0,c)$
having the Lebesgue measure $c_0 v_d c^d$. Using the previous
displays we obtain
    \begin{equation}\label{wdw}
    \min_{\|W\|=1} W^T B W \ge c\mu_{\min} \min_{\|W\|=1;S\in\calA}
        \int_{S} \big( \sum_{|s|\le \lbeta} W_s u^s \big)^2 du
        \triangleq 2\mu_0.
    \end{equation}
By the compactness argument, the minimum in (\ref{wdw}) exists and
is strictly positive.

For $i=1,\dots,n$ and any multi-indices $s_1, s_2$ such that
$|s_1|,|s_2| \le \lbeta$, define
    \begar
    T_i \triangleq \frac{1}{h^d} \xih^{s_1+s_2} \kxi
        - \int_{\dsR^d} u^{s_1+s_2} K(u) \mu(x+hu)du.
    \endar
We have $\dsE T_i = 0$,
    $|T_i| \le h^{-d}\sup_{u\in\dsR^d} \big(1+\|u\|^{2\beta}\big) K(u)
    \triangleq \kap_1h^{-d}$
and the following bound for the variance of $T_i$:
    \begar
    \Var T_i & \le & \frac{1}{h^{2d}} \dsE \xih^{2s_1+2s_2} K^2\xih\\
    & = & \frac{1}{h^{d}} \int_{\dsR^d} u^{2s_1+2s_2} K^2(u) \mu(x+hu) du\\
    & \le & \frac{\mu_{\max}}{h^d}\int_{\dsR^d} \big(1+\|u\|^{4\beta}\big) K^2(u) du
        \triangleq \frac{\kap_2}{h^d}.
    \endar
From Bernstein's inequality, we get
    \begar
    \Pn\left( |\barO_{s_1,s_2}-B_{s_1,s_2}|>\eps \right)
        = \Pn\left( \left|\frac{1}{n}\sum_{i=1}^n T_i\right|>\eps \right)
        \le 2\exp\left\{ -\frac{nh^d\eps^2}{2\kap_2+2\kap_1\eps/3} \right\}.
    \endar
This and (\ref{eq:lamo}) -- (\ref{wdw}) imply that
    \begarlab{eq:lamb}
    \Pn(\lam_{\barO} \le \mu_0) \le 2 M^2 \exp\big(-C nh^d\big)
    \endarlab
where $M^2$ is the number of elements of the matrix $\barO$. Assume
in what follows that $n$ is large enough, so that $\mu_0 > (\log
n)^{-1}$. Then for $\lam_{\barO}
> \mu_0$ we have $|\heta^*(x)-\eta(x)|\le |\heta^{LP}(x)-\eta(x)|$. Therefore,
    \begarlab{eq:etadec}
    \Pn\big(\big|\heta^*(x)-\eta(x)\big| \ge \delta \big)
        \le \Pn\big(\lam_{\barO} \le \mu_0 \big)
        + \Pn\big(\big|\heta^{LP}(x)-\eta(x)\big| \ge \delta
        , \ \lam_{\barO} > \mu_0\big).
    \endarlab
We now evaluate the second probability on the right hand side of
(\ref{eq:etadec}). For $\lam_{\barO}
> \mu_0$ we have
    $\heta^{LP}(x) = U^T(0) Q^{-1} V$
\big(where $V$ is given by \eqref{eq:v}\big). Introduce the matrix
$Z \triangleq \big(Z_{i,s}\big)_{1\le i\le n,|s|\le \lbeta}$ with
elements
    \begar
    Z_{i,s} \triangleq (X_i-x)^s \sqrt{\kxi}.
    \endar
The $s$-th column of $Z$ is denoted by $Z_s$, and we introduce
$Z^{(\eta)} \triangleq \sum_{|s|\le \lbeta} \frac{\eta^{(s)}(x)}{s!}
Z_s.$
Since $Q=Z^T Z$, we get
    \begar
    \forall |s|\le \lbeta: \quad U^T(0) Q^{-1}Z^TZ_s = \ds1_{\{s=(0,\dots,0)\}},
    \endar
hence
    $
    U^T(0) Q^{-1}Z^T Z^{(\eta)} = \eta(x).
    $
So we can write
    \begar
    \heta^{LP}(x) - \eta(x) = U^T(0) Q^{-1} (V-Z^T Z^{(\eta)})
        = U^T(0) \barO^{-1} {\bf a}
    \endar
where ${\bf a}\triangleq \frac{1}{nh^d} H(V-Z^T
Z^{(\eta)})\in\dsR^M$ and $H$ is a diagonal matrix $H \triangleq
\big(H_{s_1,s_2}\big)_{|s_1|,|s_2|\le \lbeta}$ with
    $H_{s_1,s_2} \triangleq h^{-s_1} \ds1_{\{s_1=s_2\}}.$ For $\lam_{\barO}
> \mu_0$ we get
    \begarlab{eq:etadif}
    \big| \heta^{LP}(x) - \eta(x) \big| \le \|\barO^{-1} {\bf a}\|
        \le \lam_{\barO}^{-1} \|{\bf a}\| \le \mu_0^{-1} \|{\bf a}\|
        \le \mu_0^{-1}M\max_s |a_s|,
    \endarlab
where $a_s$ are the components of the vector ${\bf a}$ given by
    \begar
    a_s = \frac{1}{nh^d} \sum_{i=1}^n \big[Y_i-\eta_x(X_i)\big] \xih^s \kxi.
    \endar
Define
    \begar
    T^{(s,1)}_i & \triangleq & \frac{1}{h^d}\big[Y_i-\eta(X_i)\big]
    \xih^s \kxi,\\
    T^{(s,2)}_i &
    \triangleq & \frac{1}{h^d}\big[\eta(X_i)-\eta_x(X_i)\big]
    \xih^s \kxi .
    \endar
We have
    \begarlab{eq:xx}
    |a_s| \le \big| \frac{1}{n}\sum_{i=1}^n T^{(s,1)}_i \big|
        + \big|\frac{1}{n}\sum_{i=1}^n \big[ T^{(s,2)}_i - \dsE T^{(s,2)}_i \big] \big|
        + \big| \dsE T^{(s,2)}_i \big|.
    \endarlab
Note that $\dsE T^{(s,1)}_i = 0$, $\big| T^{(s,1)}_i \big|  \le
\kap_1 h^{-d}$, and
    \begar
    \Var T^{(s,1)}_i & \le & 4^{-1}h^{-d}\int u^{2s} K^2(u) \mu(x+hu) du
        \le (\kap_2/4) h^{-d},\\
    \big| T^{(s,2)}_i - \dsE T^{(s,2)}_i \big|
        & \le & L \kap_1 h^{\beta-d} +  L \kap_2 h^{\beta}
        \le C h^{\beta-d},\\
    \Var T^{(s,2)}_i & \le & h^{-d}L^2\int h^{2\beta}\|u\|^{2s+2\beta}
    K^2(u) \mu(x+hu) du\le L^2\kap_2 h^{2\beta-d}.\\
    \endar
From Bernstein's inequality, for any $\eps_1,\eps_2>0$, we obtain
    \begar
    \Pn\left( \left|\frac{1}{n}\sum_{i=1}^n T_i^{(s,1)}\right| \ge \eps_1 \right)
        \le 2\exp\left\{ -\frac{nh^d\eps_1^2}{\kap_2/2+2\kap_1\eps_1/3} \right\}
    \endar
and
    \begar
    \Pn\left( \left|\frac{1}{n}\sum_{i=1}^n
        \big[ T^{(s,2)}_i - \dsE T^{(s,2)}_i \big] \right| \ge \eps_2 \right)
        \le 2\exp\left\{ -\frac{n h^d \eps_2^2}{2 L^2\kap_2 h^{2\beta}
            + 2 C h^{\beta} \eps_2/3} \right\}.
    \endar
Since also
$$
\big| \dsE T^{(s,2)}_i \big|  \le  L h^{\beta} \int \|u\|^{s+\beta}
K^2(u) \mu(x+hu) du
        \le L \kap_2 h^{\beta}
$$
we get, using (\ref{eq:xx}), that if $3\mu_0^{-1}M L\kap_2h^\beta
\le \delta\le 1$ the following inequality holds
    \begar
    \Pn\left( |a_s| \ge \frac{\mu_0\delta}{M} \right)
        & \le & \Pn\left( \left|\frac{1}{n}\sum_{i=1}^n T_i^{(s,1)}\right|
        > \frac{\mu_0\delta}{3M} \right)
        + \Pn\left( \left|\frac{1}{n}\sum_{i=1}^n
        \big[ T^{(s,2)}_i - \dsE T^{(s,2)}_i \big] \right|>\frac{\mu_0\delta}
        {3M} \right)\\
    & \le & 4 \exp\left( -C n h^d \delta^2 \right).
    \endar
Combining this inequality with \eqref{eq:lamb}, \eqref{eq:etadec} and \eqref{eq:etadif},
we obtain
    \begarlab{eq:x}
    \Pn\big(\big|\heta^*(x)-\eta(x)\big| \ge \delta \big) \le
        C_1 \exp\Big( -C_2 n h^d \delta^2 \Big)
    \endarlab
for $3m^{-1}M L\kap_2h^\beta \le \delta$ (for $\delta>1$ inequality
(\ref{eq:x}) is obvious since $\heta^*,\eta$ take values in
$[0,1]$). The constants $C_1,C_2$ in (\ref{eq:x}) do not depend on
the distribution $\px$, on its support $A$ and on the point $x\in
A$, so that we get (\ref{expon1}). Now, (\ref{expon1}) implies
(\ref{expon2}) for $Cn^{-\frac{\beta}{2\beta+d}}\le \delta$, and
thus for all $\delta
>0$ (with possibly modified constants $C_1$ and $C_2$).

\subsection{Proof of Theorems \ref{th:lb} and \ref{th:lb2}} \label{sec:prooflb}

The proof of both theorems is based on Assouad's lemma [see, e.g.,
Korostelev and Tsybakov (1993), Chapter 2 or Tsybakov (2004b),
Chapter 2]. We apply it in a form adapted for the classification
problem (Lemma 5.1 in Audibert (2004)).

For an integer $q\ge 1$ we consider the regular grid on $\dsR^d$
defined as
    $$G_q \triangleq \left\{ \left(\frac{2k_1+1}{2q},
    \dots,\frac{2k_d+1}{2q}\right)
        : k_i \in \{0,\dots,q-1\}, i=1,\dots,d \right\}.$$
Let $n_q(x) \in G_q$ be the closest point to $x\in\dsR^d$ among
points in $G_q$ (we assume uniqueness of $n_q(x)$: if there exist
several points in $G_q$ closest to $x$ we define $n_q(x)$ as the one
which is closest to 0). Consider the partition
$\calX'_1,\dots,\calX'_{q^d}$ of $[0,1]^d$ canonically defined using
the grid $G_q$ ($x$ and $y$ belong to the same subset if and only if
$n_q(x)=n_q(y)$). Fix an integer $m \leq q^d$. For any
$i\in\{1,\dots,m\}$, we define $\calX_i \triangleq \calX'_i$ and
$\calX_0 \triangleq \dsR^d\setminus\cup_{i=1}^m \calX_i$, so that
$\calX_0,\dots,\calX_{m}$ form a partition of~$\dsR^d$.

Let $u:\dsR_+\rightarrow\dsR_+$ be a nonincreasing infinitely differentiable function such that
    $u=1$ on $[0, 1/4]$ and $u=0$ on $[1/2, \infty)$.
    One can take, for example,
$u(x) = \Big(\int_{1/4}^{1/2} u_1(t)dt\Big)^{-1} \int_x^{\infty}
u_1(t)dt$ where the infinitely differentiable function $u_1$ is
defined as
    \begar
    u_1(x)=\left\{ \begin{array}{ll}
        \exp\Big\{ -\frac{1}{(1/2-x)(x-1/4)} \Big\} & \qquad \text{for }
        x\in (1/4,1/2),\\
        0 & \qquad \text{otherwise.}
        \end{array} \right.
    \endar
Let $\phi:\dsR^d\rightarrow\dsR_+$ be the function defined as
    \begar
    \phi(x) \triangleq C_\phi u(\|x\|),
    \endar
where the positive constant $C_\phi$ is taken small enough so ensure
that
    $|\phi(x') - \phi_x(x')| \le L \|x'-x\|^{\beta}$
    for any $x,x'\in\dsR^d$.
Thus, $\phi\in \Sigma( \beta , L , \dsR^d )$.

Define the hypercube ${\cal H}=\big\{ \dsP_{\sigmav} :
\sigmav=(\sigma_1,\dots,\sigma_m) \in \{-1,1\}^m \big\}$ of
probability distributions $\dsP_{\sigmav}$ of $(X,Y)$ on $\calZ =
\dsR^d\times\{0,1\}$ as follows.

For any $\dsP_{\sigmav}\in{\cal H}$ the marginal distribution of $X$
does not depend on $\sigmav$, and has a density $\mu$ w.r.t. the
Lebesgue measure on $\dsR^d$ defined in the following way. Fix
$0<w\le m^{-1}$ and a set $A_0$ of positive Lebesgue measure
included in $\calX_0$ (the particular choices of $A_0$ will be
indicated later), and take: (i) $\mu(x)=w/\lam[\calB(0,(4q)^{-1})]$
if $x$ belongs to a ball $\calB(z,(4q)^{-1})$ for some $z\in G_d$,
(ii) $\mu(x)=(1-mw)/\lam[A_0]$ for $x\in A_0$, (iii) $\mu(x)=0$ for
all other $x$.

Next, the distribution of $Y$ given $X$ for $\dsP_{\sigmav}\in{\cal
H}$ is determined by the regression function
$\eta_{\sigmav}(x)=P(Y=1|X=x)$ that we define as $\eta_{\sigmav}(x)
= \frac{1+\sigma_j \varphi(x)}{2}$ for any $x\in\calX_j$,
$j=1,\dots,m$, and $\eta_{\sigmav} \equiv 1/2$ on $\calX_0$, where
    $\varphi(x) \triangleq q^{-\beta} \phi\big( q [x-n_q(x)] \big).$
We will assume that $C_\phi \le 1$ to ensure that $\varphi$ and
$\eta_{\sigmav}$ take values in $[0,1]$.

For any $s\in\dsN^d$ such that $|s|\le\fbeta$, the partial
derivative $D^s\varphi$ exists, and
    $D^s\varphi(x) = q^{|s|-\beta} D^s\phi\big( q [x-n_q(x)] \big)$.
Therefore, for any $i\in\{1,\dots,m\}$ and any $x,x'\in\calX_i$, we
have
        \begin{equation*}
        |\varphi(x') - \varphi_x(x')| \leq L \|x-x'\|^{\beta}.
        \end{equation*}
This implies that for any $\sigmav \in \{-1,1\}^m$ the function
$\eta_{\sigmav}$ belongs to the H\"{o}lder class
$\Sigma\big(\beta,L,\dsR^d\big)$.

We now check the margin assumption. Set
$x_0=\big(\frac{1}{2q},\dots,\frac{1}{2q}\big)$. For any $\sigmav
\in \{-1,1\}^m$ we have
    \begar
    \dsP_{\sigmav}\big( 0<\big| \eta_{\sigmav}(X) - 1/2 \big| \leq t
    \big)
        & = & m \dsP_{\sigmav}\big(0<\phi[q(X-x_0)] \leq 2t q^\beta \big)\\
    & = & m \int_{\calB(x_0,(4q)^{-1})}
    \ds1_{\{0<\phi[q(x-x_0)] \leq 2t q^\beta\}}
    \frac{w}{\lam[\calB(0,(4q)^{-1})]} dx\\
    & = & \frac{mw}{\lam[\calB(0,1/4)]} \int_{\calB(0,1/4)}
    \ds1_{\{\phi(x) \leq 2tq^\beta\}} dx\\
    & = & mw \ds1_{\{t \geq C_\phi/(2q^\beta)\}}.
    \endar
Therefore, the margin assumption (MA) is satisfied if $mw =
O(q^{-\alpha \beta}).$

According to Lemma 5.1 in Audibert (2004), for any classifier $\hf$
we have
        \begin{equation}\label{lba}
        \underset{P \in {\cal H}}{\sup}
        \big\{ \dsE R( \hf ) - R( f^* ) \big\}
            \geq mwb' (1-b\sqrt{nw})/2
        \end{equation}
where \begar b & \triangleq & \left[1-\big(\int_{\calX_1}
\sqrt{1-\varphi^2(x)} \, \mu_1(x)dx\big)^2\right]^{1/2}
        = C_\phi q^{-\beta},\\
    b' & \triangleq & \int_{\calX_1} \varphi(x)\mu_1(x)dx
    =
    C_\phi q^{-\beta}
    \endar
with $\mu_1(x) = \mu(x)/\int_{\calX_1} \mu(z)dz$.

We now prove Theorem \ref{th:lb}. Take $q=\big\lfloor \bar{C}
n^{\frac{1}{2\beta+d}} \big\rfloor$, $w=C'q^{-d}$ and $m=\big\lfloor
C''q^{d-\alpha \beta} \big\rfloor$ with some positive constants
$\bar{C}, C', C''$ to be chosen, and set
$A_0=[0,1]^d\setminus\cup_{i=1}^m \calX_i$.  The condition $\alpha
\beta \le d$ ensures that the above choice of $m$ is not degenerate:
we have $m\ge 1$ for $C''$ large enough. We now prove that ${\cal
H}\subset \calPS$ under the appropriate choice of $\bar{C},C',C''$.
In fact, select these constants so that the triplet $(q,w,m)$ meets
the conditions $m\le q^d$, $0<w\le m^{-1}$, $mw = O(q^{-\alpha
\beta})$. Then, in view of the argument preceding (\ref{lba}), for
any $\sigmav \in \{-1,1\}^m$ the regression function
$\eta_{\sigmav}$ belongs to $\Sigma\big(\beta,L,\dsR^d\big)$ and
Assumption (MA) is satisfied. We now check that $\px$ obeys the
strong density assumption. First, the density $\mu(x)$ equals to a
positive constant for $x$ belonging to the union of balls
$\cup_{i=1}^m\calB(z_i,(4q)^{-1})$ where $z_i$ is the center of
$\calX_i$, and $\mu(x)=(1-mw)/(1-mq^{-d})=1+o(1)$, as $n\to\infty$,
for $x\in A_0$. Thus, $\mu_{\min}\le \mu(x)\le \mu_{\max}$ for some
positive $\mu_{\min}$ and $\mu_{\max}$. (Note that this construction
does not allow to choose any prescribed values of $\mu_{\min}$ and
$\mu_{\max}$, because $\mu(x)=1+o(1)$. The problem can be fixed via
a straightforward but cumbersome modification of the definition of
$A_0$ that we skip here.) Second, the $(c_0,r_0)$-regularity of the
support $A$ of $\px$ with some $c_0>0$ and $r_0>0$ follows from the
fact that, by construction, $\lam (A\cap \calB(x,r)) = (1+o(1))\lam
([0,1]^d\cap \calB(x,r))$ for all $x\in A$ and $r>0$ (here again we
skip the obvious generalization allowing to get any prescribed
$c_0>0$). Thus, the strong density assumption is satisfied, and we
conclude that ${\cal H}\subset \calPS$. Theorem \ref{th:lb} now
follows from (\ref{lba}) if we choose $C'$ small enough.

Finally, we prove Theorem \ref{th:lb2}. Take $q=\big\lfloor C
n^{\frac{1}{(2+\alpha)\beta+d}} \big\rfloor$, $w= C' q^{2\beta}/n$
and $m=q^d$ for some constants $C>0$, $C'>0$, and choose $A_0$ as a
Euclidean ball contained in $\calX_0$. As in the proof of Theorem
\ref{th:lb}, under the appropriate choice of $C$ and $C'$, the
regression function $\eta_{\sigmav}$ belongs to
$\Sigma\big(\beta,L,\dsR^d\big)$ and the margin assumption (MA) is
satisfied. Moreover, it is easy to see that the marginal
distribution of $X$ obeys the mild density assumption (the
$(c_0,r_0)$-regularity of the support of $\px$ follows from
considerations analogous to those in the proof of Theorem
\ref{th:lb}). Thus, ${\cal H}\subset \calPS'$. Choosing $C'$ small
enough and using (\ref{lba}) we obtain Theorem \ref{th:lb2}.

\subsection{Proof of Proposition \ref{th:hitcross}} \label{sec:proofhc}

The following lemma describes how the smoothness constraint on the
regression function $\eta$ at some point $x\in\dsR^d$ implies that
$\eta$ ``stays close'' to $\eta(x)$ in the vicinity of $x$.
\begin{lem} \label{lem:hitcross}
For any distribution $P\in \calPS$ with regression function $\eta$
and for any $\kap>0$, there exist $L'>0$ and $t_0>0$ such that for
any $x$ in the support of $\px$ and $0<t\le t_0$, we have
    $$\px\Big[ \big| \eta(X)-\eta(x)\big| \le t ;
    X\in\calB\big(x,\kap t^{\frac{1}{1\wedge \beta}} \big)\Big] \ge L' t^{\frac{d}{1\wedge \beta}}.$$
\end{lem}

\begin{proof}[Proof of Lemma \ref{lem:hitcross}]
Let $A$ denote the support of $\px$.
Let us first consider the case $\beta\le 1$. Then for any $x,x'\in\dsR^d$, we have
    $\big| \eta(x') - \eta(x) \big| \le L \| x'-x \|^\beta.$
Let $\kap'= \kap \wedge L^{-1/\beta}$.
For any $0<t\le L r_0^\beta$, we get
    \begar
    \px\Big[ \big| \eta(X) - \eta(x) \big| \le t ; X\in\calB\big(x,\kap t^{\frac{1}{1\wedge \beta}} \big) \Big]\\
    \qquad\qquad =
        \px\Big[ \big| \eta(X) - \eta(x) \big| \le t ;
        X \in \calB\big(x,\kap t^{\frac{1}{\beta}} \big) \cap A \Big]\\
    \qquad\qquad  \ge
        \px\Big[ X \in \calB\big(x,\kap t^{\frac{1}{\beta}} \wedge \big(\frac{t}{L}\big)^{\frac{1}{\beta}} \big)
        \cap A \Big]\\
    \qquad\qquad  \ge
        \mu_{\min} \lam\Big[ \calB\big(x, \kap' t^{\frac{1}{\beta}}  \big) \cap A \Big]\\
    \qquad\qquad  \ge
        c_0 \mu_{\min} \lam\Big[ \calB\big(x, \kap' t^{\frac{1}{\beta}} \big) \Big]\\
    \qquad\qquad  \ge
        c_0 \mu_{\min} v_d (\kap')^d t^{\frac{d}{\beta}},
    \endar
which is the desired result with $L' \le c_0 \mu_{\min} v_d (\kap')^d$ and $t_0\le L r_0^\beta$.

For the case $\beta>1$, by assumption, $\eta$ is continuously
differentiable. Let $\calC(A)$ be the convex hull of the support $A$
of $\px$. By compactness of $\calC(A)$, there exists $C>0$ such that
for any $s\in\dsN^d$ with $|s|=1$,
    $$\underset{x\in\calC(A)}{\sup} \; \big|D^s \eta(x) \big| \le C.$$
So we have for any $x,x' \in A$,
    $$|\eta(x)-\eta(x')| \le C \|x-x'\|.$$
The rest of the proof is then similar to the one of the first case.
\end{proof}

\begin{itemize}
\item We will now prove the first item of Proposition \ref{th:hitcross}.
Let $P\in\calPS$ such that the regression function associated with
$P$ hits $1/2$ at $x_0\in \overset{\circ}{A}$, where $\overset{\circ}{A}$ denotes the
interior of the support of $\px$.
Let $r>0$ such that $\calB(x_0,r) \subset A$. Let $x\in \calB(x_0,r)$ such that
$\eta(x)\neq \demi$. Let $t_1=\big|\eta(x)-1/2\big|.$
For any $0< t\le t_1$, let $x_t\in[x_0;x]$ such that $\big| \eta(x_t) - 1/2 \big| = t/2.$
We have $x_t\in A$ so that we can apply Lemma \ref{lem:hitcross} (with $\kappa=1$
for instance) and obtain for any $0< t\le t_1 \wedge (4t_0)$
    $$\px\Big[ 0 < \big| \eta(X) - 1/2 \big| \le t \Big]
        \ge \px\Big[ \big| \eta(X) - \eta(x_t) \big| \le t/4 \Big] \ge L' (t/4)^{\frac{d}{1\wedge \beta}}.$$
Now from the margin assumption, we get that for any small enough $t>0$
    $C_0 t^\alpha \ge L' (t/4)^{\frac{d}{1\wedge \beta}},$
hence $\alpha \le \frac{d}{1\wedge \beta}.$

\item For the second item of Proposition \ref{th:hitcross}, to skip cumbersome details,
we may assume that $\calC$ contains the unit ball in $\dsR^d$.
Consider the distribution such that
    \begin{itemize}
    \item $\px$ is the uniform measure on
    $\big\{ (x_1,\dots,x_d) \in\dsR^d : |x_1-1/4| + |x_2| + \cdots + |x_d| \le 1/4 \big\}$
    \item the regression function associated with $P$ is
    $$\eta(x_1,\dots,x_d)=\frac{1+C_\eta \sign(x_1) |x_1|^{\beta\wedge 1} u(x_1)}{2},$$ where
    $$u(t)=\left\{ \begin{array}{lll}
    \exp\big(-\frac{1}{1-t^2}\big) & \text{if } |t|<1\\
    0 & \text{otherwise},
    \end{array} \right.$$
    and $0<C_\eta\le 1$ is small enough so that for any $x,x'\in\dsR^d$, $\eta$ satisfies
    $$|\eta(x')-\eta_x(x')| \le L \|x-x'\|^\beta.$$
    \end{itemize}
For appropriate positive parameters $c_0,r_0,\mu_{\max}>\mu_{\min}>0$,
the only non-trivial task in checking that $P$ belongs to $\calPS$ is to
check the margin assumption. For $t$ small enough, we have
    $$\px\Big[ \big| \eta(X)-1/2\big| \le t \Big]
        \le \px\Big[ | X_1 |^{\beta \wedge 1} \le C t ;
        |X_1-1/4| + |X_2| + \cdots + |X_d| \le 1/4 \Big]$$
for some $C>0$. Therefore, we have
    $\px\Big[ 0 < \big| \eta(X)-1/2\big| \le t \Big] \le C t^{\frac{d}{\beta\wedge 1}}$.
So the margin assumption is satisfied for an appropriate $C_0$ whenever $\alpha \le \frac{d}{\beta \wedge 1}.$
Since $\eta$ hits $1/2$ at $0_{\dsR^d}$ which is in boundary of the support of $\px$, we have proved the second assertion.

\item For the third assertion of Proposition \ref{th:hitcross},
to avoid cumbersome details again, we may assume that $\calC$ contains the unit ball in $\dsR^d$.
Consider the distribution such that
    \begin{itemize}
    \item $\px$ is the uniform measure on the unit ball,
    \item the regression function associated with $P$ is
    $$\eta(x)=\frac{1+C_\eta \|x\|^2 u(\|x\|^2/2)}{2},$$ where
    $0<C_\eta\le 1$ is small enough so that for any $x,x'\in\dsR^d$, $\eta$ satisfies
    $$|\eta(x')-\eta_x(x')| \le L \|x-x'\|^\beta.$$
    \end{itemize}
For appropriate positive parameters $C_0,c_0,r_0,\mu_{\max}>\mu_{\min}>0$,
the distribution $P$ belongs to $\calPS$ provided that $\alpha\le d/2$ (in order that
the margin assumption holds). We have obtained the desired result
since $\eta$ hits $1/2$ at $0_{\dsR^d}$ which is in the interior of the support of $\px$.

\item For the last item of Proposition \ref{th:hitcross}, let $P\in \calPS$
such that the regression function $\eta$ associated with $P$ crosses
$1/2$ at $x_0\in \overset{\circ}{A}$. For $d=1$, from
the first item of the theorem, we necessarily have $\alpha(\beta
\wedge 1) \le 1$. Let us now consider the case: $d>1$.

Figure \ref{fig:hit} will help to keep track of the
following notation.

\begin{figure}[h]
\begin{center}
\includegraphics{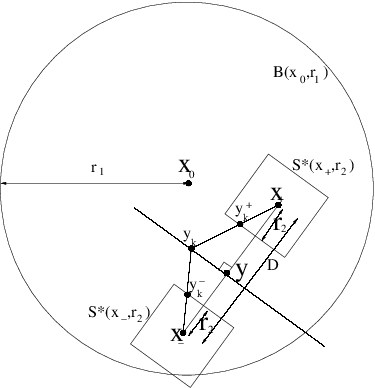}
\caption{Notation summary}
\label{fig:hit}
\end{center}
\end{figure}

Let $r_1>0$ such that $\calB(x_0,3r_1) \subset A$.
Introduce $x_-$ and $x_+$ in $\calB(x_0,r_1)$ such that
    $\eta(x_-)<1/2$ and $\eta(x_+)>1/2.$
Let $t_1= \big( 1/2-\eta(x_-)\big)\wedge\big( \eta(x_+)- 1/2\big).$
Define $y=\frac{x_-+x_+}{2}$, $e_d=\frac{x_+-x_-}{\|x_+-x_-\|}$
and $D=\|x_+-x_-\|$. Let $e_1,\dots,e_{d-1}$ be unit vectors such that
$e_1,\dots,e_{d}$ is an orthonormal basis of $\dsR^d$. Let $\calB^*(x,r)$
(resp. $\calS^*(x,r)$) denote the ball (resp. the sphere) centered at $x$ and of radius
$r$ wrt the norm $\|x\|_* = \sup_{1 \le i \le d} |\langle x , e_i \rangle|.$

Since $\eta$ is continuous, there exists $r_2>0$ such that
    $$\left\{ \begin{array}{lll}
    \eta(x) < 1/2 - t_1/2 & \text{for any } x\in \calB^*(x_-,r_2)\\
    \eta(x) > 1/2 + t_1/2& \text{for any } x\in \calB^*(x_+,r_2)
    \end{array} \right.$$
Let $\zeta=\frac{1}{\beta \wedge 1}$. For any $k=(k_1,\dots,k_{d-1}) \in \dsZ^{d-1}$, introduce
    $$y_k = y + t^\zeta \sum_{i=1}^{d-1} k_i e_i.$$
For any $t$ in $]0;t_1[$, consider the grid
    $G = \big\{ y_k ; k\in\dsZ^{d-1},
    \underset{1 \le i \le d-1}{\max} |k_i| \le \frac{D}{2\sqrt{d-1}t^{\zeta}} \big\}.$
For any $y_k$ in $G$, we have $\|y_k-y\| \le \sqrt{d-1} \underset{1 \le i \le d-1}{\max} |t^\zeta k_i| \le D/2 \le r_1$.
Therefore, using that $y\in\calB(x_0,r_1)$, the grid $G$ is included in $\calB(x_0,2 r_1)$.
For any $y_k\in G$, let $y_k^-=[x_-;y_k] \cap \calS^*(x_-,r_2)$ and
$y_k^+=[x_+;y_k] \cap \calS^*(x_+,r_2)$. Since $\|y_k-y\| \le D/2$,
we have
    $y_k^- = x_- + r_2 e_d + \frac{2 r_2}{D} t^\zeta \sum_{i=1}^{d-1} k_i e_i$
and
    $y_k^+ = x_+ - r_2 e_d + \frac{2 r_2}{D} t^\zeta \sum_{i=1}^{d-1} k_i e_i$.

For any $y_k$ in $G$, consider the continuous path formed by the segments
$[y_k^-;y_k]$ and $[y_k;y_k^+]$. Since $\eta$ is continuous on this path,
there exists $w_k \in \gamma_k\triangleq [y_k^-;y_k] \cup [y_k;y_k^+]$ such that
    $\eta(w_k)=1/2+t/2$.
Now let us show that when $k\neq k'$, $w_k$ and $w_{k'}$ are at least $\frac{\sqrt{2} r_2 }{D} t^\zeta$ away from each other.
The distance between $w_k$ and $w_{k'}$ is not less than the distance between the paths $\gamma_k$ and
$\gamma_{k'}$. Let $U$ denote the biggest integer smaller than or equal to $\frac{D}{2\sqrt{d-1}t^{\zeta}}$.
When $y_k\neq y_{k'}$ in $G$, the distance between $\gamma_k$ and $\gamma_{k'}$ is minimum for
    $k=K \triangleq (U,\dots,U)$ and $k'= K' \triangleq (U-1,U,\dots,U)$.
This distance is equal to the distance between $y_K^-$ and its orthogonal projection on $[y_{K'}^-;y_{K'}]$,
which is the distance between $y_K^-$ and the line $(x_-;y_{K'})$. Let $K''=(0,U,\dots,U)\in \dsZ^{d-1}$.
To compute this distance $V$, it suffices to look at the plane $(x_-;y_{K''};y_K)$ (see figure \ref{fig:hit2}).

\begin{figure}[h]
\begin{center}
\includegraphics{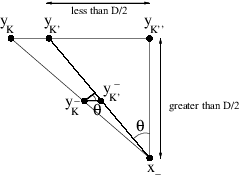}
\caption{plane $(x_-;y_{K''};y_K)$}
\label{fig:hit2}
\end{center}
\end{figure}

We obtain that the angle $\theta$ between $y_{K'}-x_-$ and $y_{K''}-x_+$ is smaller than $\pi/4$.
As a consequence, $V= \|y_{K}^- - y_{K'}^-\| \cos \theta \ge \sqrt{2} r_2 t^\zeta / D$.

Finally, focusing on the behaviour of the regression function near the $w_k$'s,
by using Lemma \ref{lem:hitcross} with $\kappa=\frac{4^\zeta r_2}{\sqrt{2} D}$,
we obtain that there exists $L'>0$ and $t_0>0$ such that
for any $0<t<4t_0 \wedge t_1$,
    $$\begin{array}{lll}
    C_0 t^\alpha & \ge & \px\Big[ 0<\big| \eta(X) - \demi \big| \le t \Big]\\
    & \ge & \underset{k\in\dsZ^{d-1}: \underset{1 \le i \le d-1}{\max} |k_i| \le \frac{D}{2\sqrt{d-1}t^{\zeta}}}{\sum}
        \px\Big[ \big| \eta(X) - \eta(w_k) \big| \le t/4 ;
        X \in \calB\big(w_k,\frac{r_2 t^\zeta}{\sqrt{2}D}\big)\Big]\\
    & \ge & (2U+1)^{d-1} L' (t/4)^{d \zeta}\\
    & \ge & \big(\frac{D}{2\sqrt{d-1}t^{\zeta}}\big)^{d-1} L' (t/4)^{d \zeta}\\
    & \ge & C t^{\zeta},
    \endar
hence $\alpha \le \zeta$ (which is the desired result).

For the converse, the proof is similar to the ones of the second and third assertions of
the proposition. Without loss of generality, we may assume that
$\calS=\big\{ (x_1,\dots,x_d) \in\dsR^d : \underset{1 \le i \le d}{\max} |x_i| \le 1/2 \big\}$
is a subset of $\calC$.
we consider the distribution $P$ such that
    \begin{itemize}
    \item $\px$ is the uniform measure on $\calS$
    \item the regression function associated with $P$ is
    $$\eta(x_1,\dots,x_d)=\frac{1+C_\eta \sign(x_1) |x_1|^{\beta\wedge 1} u(x_1)}{2},$$ where
    $0<C_\eta \le 1$ is small enough so that for any $x,x'\in\dsR^d$, $\eta$ satisfies
    $$|\eta(x')-\eta_x(x')| \le L \|x-x'\|^\beta.$$
    \end{itemize}
For small enough $t>0$, we have
    $$\px\Big[ \big| \eta(X)-1/2\big| \le t \Big]
        \le \px\big[ | X_1 |^{\beta \wedge 1} \le C t \big],$$
for some constant $C>0$, so that we have
    $\px\Big[ 0 < \big| \eta(X)-\demi\big| \le t \Big] \le 2 (C t)^{\frac{1}{\beta\wedge 1}}.$
As a consequence, for appropriate parameters $C_0,c_0,r_0,\mu_{\max}>\mu_{\min}>0$,
the distribution $P$ belongs to $\calPS$ whenever $\alpha \le
\frac{1}{\beta \wedge 1}.$ Since $\eta$ crosses $1/2$ at
$0_{\dsR^d}$ which is in the interior of the support of $\px$, the
converse holds.

\end{itemize}

\subsection{Proof of Theorem \ref{th:ub3}} \label{proof:ub2}

We prove the theorem for $p<\infty$. The proof for $p=\infty$ is
analogous. For any decision rule $f$ we set $d(f) \triangleq
R(f)-R(f^*)$ and
$$
f^{**}(x,f)\triangleq\left\{
\begin{array}{lcr}
 f^*(x) & \text{if}
&\eta(x)\neq 1/2,\\
f(x) & \text{if} &\eta(x)= 1/2,
\end{array}\right. \qquad \forall \ x\in\dsR^d.
$$
\begin{lem}\label{mts}
Under Assumption (MA) for any decision rule $f$ we have
\begin{equation}\label{p1t}
\px(f(X)\neq f^{**}(X,f))\le C d(f)^{\alpha/(1+\alpha)}.
\end{equation}
\end{lem}
\begin{proof} Note that $f^{**}(\cdot,f)$ is a Bayes rule,
and following the same lines as in Proposition 1 of Tsybakov (2004a)
we get $\px(f(X)\neq f^{**}(X,f),\, \eta(X)\neq 1/2)\le C
d(f)^{\alpha/(1+\alpha)}$. It remains to observe that $\px(f(X)\neq
f^{**}(X,f), \,\eta(X)\neq 1/2) = \px(f(X)\neq f^{**}(X,f))$.
\end{proof}

For a Borel function $\bar\eta$ on $\dsR^d$ define $
f_{\bar\eta}\triangleq\ds1_{\{\bar\eta \geq 1/2\}}$,
$f^*_{\bar\eta}(\cdot)\triangleq f^{**}(\cdot,f_{\bar\eta})$ and
$$
Z_n(f_{\bar\eta}) \triangleq [R_n(f_{\bar\eta})-R_n(f^*_{\bar\eta})]
- [R(f_{\bar\eta})-R(f^*_{\bar\eta})] =
[R_n(f_{\bar\eta})-R_n(f^*_{\bar\eta})] - d(f_{\bar\eta}).
$$

Let $\eta_n$ be an element of ${\calN}_{\ve_n}$ such that
$\|\eta_n-\eta\|_{p,\lam} \le \ve_n$, where $\|\cdot\|_{p,\lam}$ is
the $L_p(\calC,\lam)$ norm. In view of the assumption on $\calP$ we
have $\|\eta_n-\eta\|_{p} \le \mu_{\max}^{1/p}\ve_n$ where
$\|\cdot\|_{p}$ is the $L_p(\dsR^d,\px)$ norm. It follows from the
comparison inequality (\ref{eq:comp3}) that $d(f_{\eta_n}) \le
C\ve_n^{\frac{(1+\alpha)p}{p+\alpha}} \triangleq \delta_n$. Set
$$
\Delta_n = C n^{- \frac{(1+\alpha)p}{(2+\alpha)p+\rho(p+\alpha)}}
$$
(i.e., $\Delta_n$ is of the order of desired rate). Fix $t>0$ and
introduce the set
$$
{\calN}^*_n=\{\bar\eta\in{\calN}_{\ve_n}: \ d(f_{\bar\eta})\ge
t\Delta_n\}.
$$
For any $t>0$ we have
\begin{eqnarray*}
\dsP (d(\hat f^s_n) \ge t\Delta_n)  &\le &\dsP
(\min_{\bar\eta\in{\calN}^*_n}
[R_n(f_{\bar\eta}) - R_n(f_{\eta_n})] \le 0) \\
&=&\dsP (\min_{\bar\eta\in{\calN}^*_n}
[Z_n(f_{\bar\eta}) - Z_n(f_{\eta_n}) + d(f_{\bar\eta}) - d(f_{\eta_n})]
\le 0) \\
&\le&\dsP (\min_{\bar\eta\in{\calN}^*_n} [Z_n(f_{\bar\eta}) -
Z_n(f_{\eta_n}) + d(f_{\bar\eta})/2+t\Delta_n/2 - d(f_{\eta_n})]
\le 0) \\
&\le&\dsP (\min_{\bar\eta\in{\calN}^*_n} [Z_n(f_{\bar\eta})
+d(f_{\bar\eta})/2] \le 0)\\
&& +\dsP ( Z_n(f_{\eta_n}) \ge t\Delta_n/2 - d(f_{\eta_n}))
\\
&\le&\dsP (\min_{\bar\eta\in{\calN}^*_n} [Z_n(f_{\bar\eta})
+d(f_{\bar\eta})/2] \le 0)\\
&& +\dsP ( Z_n(f_{\eta_n}) \ge t\Delta_n/2 - \delta_n).
\end{eqnarray*}
Since $\Delta_n$ is of the same order as $\delta_n$, we can choose
$t$ large enough to have $t\Delta_n/2 - \delta_n\ge t\Delta_n/4$.
Thus,
\begin{eqnarray*}
\dsP (d(\hat f^s_n) \ge t\Delta_n)  &\le & {\rm card}\ {\calN}^*_n \
\max_{\bar\eta\in{\calN}^*_n}\dsP ( Z_n(f_{\bar\eta})
 \le - d(f_{\bar\eta})/2)\\
&& +\dsP ( Z_n(f_{\eta_n}) \ge t\Delta_n/4 )\\
&\le & \exp(A'\ve_n^{-\rho})  \max_{\bar\eta\in{\calN}^*_n}\dsP (
Z_n(f_{\bar\eta})
 \le - d(f_{\bar\eta})/2)\\
&& +\dsP ( Z_n(f_{\eta_n}) \ge t\Delta_n/4 ).
\end{eqnarray*}
Note that for any decision rule $f$ the value $Z_n(f)$ is an average
of $n$ i.i.d. bounded and centered random variables whose variance
does not exceed $\px(f(X)\neq f^{**}(X,f))$. Thus, using Bernstein's
inequality and (\ref{p1t}) we obtain
$$
\dsP (-Z_n(f) \ge a) \le \exp\left(- \frac{Cna^2}{a +
d(f)^{\alpha/(1+\alpha)}}\right), \quad \forall \ a>0.
$$
Therefore, for $\bar\eta\in{\calN}^*_n$,
\begin{eqnarray*}
\dsP ( Z_n(f_{\bar\eta})
 \le - d(f_{\bar\eta})/2)&\le&
 \exp (- Cn d(f_{\bar\eta})^{(2+\alpha)/(1+\alpha)})\\
&\le&
 \exp (- Cn
 (t\Delta_n)^{(2+\alpha)/(1+\alpha)}).
\end{eqnarray*}
Similarly, for $t>C$,
\begin{eqnarray*}
\dsP ( Z_n(f_{\eta_n}) \ge t\Delta_n/4 )&\le& \exp \left(- \frac{Cn
 \Delta_n^2}{\Delta_n
 +d(f_{\eta_n})^{\alpha/(1+\alpha)}}\right)\\
 &\le& \exp \left(- \frac{Cn
 \Delta_n^2}{\Delta_n
 +\delta_n^{\alpha/(1+\alpha)}}\right)\\
  &\le& \exp \left(- Cn
 \Delta_n^{(2+\alpha)/(1+\alpha)}\right).
\end{eqnarray*}
The result of the theorem follows now from the above inequalities
and the relation $n
 \Delta_n^{(2+\alpha)/(1+\alpha)} \asymp \ve_n^{-\rho}$.

\footnotesize{

}
\begin{flushright}
\footnotesize $^1${\sc Center for Education and Research in Informatics\\}
{\sc Ecole Nationale des Ponts et Chauss\'{e}es\\}
{\sc 19, rue Alfred Nobel\\}
{\sc Cit\'e Descartes, Champs-sur-Marne\\}
{\sc 77455 Marne-La-Vallée, France\\}
{\sc e-mail}: audibert@certis.enpc.fr
\bigskip\\
$^2${\sc Laboratoire de Probabilit\'es et Mod\`eles }\\
{\sc Al\'eatoires (UMR CNRS 7599),}\\
{\sc Universit\'e Paris VI}\\
{\sc 4, pl.Jussieu, Bo\^{\i}te courrier 188,}\\
{\sc 75252 Paris, France}\\
{\sc e-mail}: tsybakov@ccr.jussieu.fr

\end{flushright}


\begin{thebibliography}{99}


\bibitem{Aud03d}
Audibert, J.-Y. (2004). {Classification using $\text{Gibbs}$
estimators under complexity and margin assumptions}. {Preprint},
{Laboratoire de Probabilit{\'e}s et Model{\`e}s Al{\'e}atoires},
{http://www.proba.jussieu.fr/mathdoc/textes/PMA-908.pdf}.

\bibitem{BJM03} Bartlett, P.L., Jordan, M.I. and McAuliffe, J.D. (2003). Convexity,
classification and risk bounds. Techn. Report 638, University of
California at Berkeley.

\bibitem{BS} Birman, M.S. and Solomjak, M.Z. (1967)
Piecewise-polynomial approximation of functions of the classes
$W^{\alpha}_p$. {\it Mat. Sbornik} {\bf 73} 295-317.

\bibitem{BBM04} Blanchard, G., Bousquet,O. and Massart, P. (2004)
Statistical performance of support vector machines. Manuscript.
http://www.kyb.mpg.de/publications/pss/ps2731.ps

\bibitem{BLV03} Blanchard, G., Lugosi, G. and Vayatis, N. (2003). On the rate of
convergence of regularized boosting classifiers. {\it Journal of
Machine Learning Research} {\bf 4} 861-894.

\bibitem{DGL96} Devroye, L. , Gy\"orfi, L. and Lugosi, G. (1996). {\it A
Probabilistic Theory of Pattern Recognition.} Springer, New York,
Berlin, Heidelberg.

\bibitem{gea02}
{Gy\"{o}rfi, L., Kohler, M., Krzy\.{z}ak, A. and Walk, H.} (2002).
{\it A Distribution-Free Theory of Nonparametric Regression}.
Springer, New York.


\bibitem{kt61}
Kolmogorov, A.N. and Tikhomorov V.M. (1961) $\epsilon$-entropy and
$\epsilon$-capacity of sets in function spaces. {\it Translations of
the American Mathematical Society} {\bf 17} 277-364.

\bibitem{K05} Koltchinskii, V. (2005) Local Rademacher complexities and oracle
inequalities in risk minimization. Manuscript.

\bibitem{KB05} Koltchinskii, V. and Beznosova, O. (2005) Exponential
convergence rates in classification. Proceedings of the 18th
Conference on Learning Theory (COLT-2005).

\bibitem{KT93} Korostelev, A. P. and Tsybakov, A. B. (1993). {\it Minimax Theory
of Image Reconstruction.} Lecture Notes in Statistics {\bf 82},
Springer, New York, Berlin, Heidelberg.

\bibitem{LW04} Lugosi, G. and Wegkamp, M. (2004). Complexity regularization via
localized random penalties. {\it Ann.\ Statist.} {\bf 32} 1679 -
1697.

\bibitem{Mam99} Mammen, E. and Tsybakov, A. B. (1999). Smooth discrimination
analysis. {\it Ann. Statist.} {\bf 27} 1808 - 1829.

\bibitem{MN03} Massart, P. and N\'ed\'elec, E. (2003) Risk bounds for
statistical learning. Preprint.
http://www.math.u-psud.fr/~massart/margin.pdf.

\bibitem{ScSt03}
Scovel, J.C. and Steinwart, I. (2004). Fast Rates for Support Vector
Machines. Los Alamos National Laboratory Technical Report LA-UR
03-9117.
http://www.c3.lanl.gov/~ingo/publications/ann-04a.pdf.

\bibitem{Sto80} Stone, C.J. (1980)
    {Optimal rates of convergence for nonparametric estimators},
    {\it Ann. \ Statist.}
    {\bf 8}
    {1348--1360}.

\bibitem{Sto82}
Stone, C.J. (1982) Optimal global rates of convergence for
nonparametric regression. {\it Ann. \ Statist.} {\bf 10}
1040-1053.


\bibitem{TaVa04}
Tarigan, B. and van de Geer, S. (2004). Adaptivity of Support Vector
Machines with $\ell_1$ Penalty. Preprint, University of Leiden,
~http://www.math.leidenuniv.nl/~geer/svm4.pdf.

\bibitem{Tsy04} Tsybakov, A.B. (2004a). Optimal aggregation of classifiers in
statistical learning. {\it Ann. Statist.} {\bf 32} 135-166.

\bibitem{tsy04}
{Tsybakov, A. B.} (2004b). {\it Introduction \`{a} l'estimation
non--param\'{e}trique.} Springer, Berlin.

\bibitem{tsvdg05}
Tsybakov, A. B. and van de Geer, S. (2005) Square root penalty:
adaptation to the margin in classification and in edge estimation.
{\it Annals of Statistics} {\bf 33} 1203 - 1224.

\bibitem{vdg00}
{van de Geer, S.} (2000). {\it Empirical Processes in M-Estimation},
Cambridge Univ. Press.

\bibitem{v2}
Vapnik, V. N. (1998). {\it Statistical Learning Theory.} Wiley, New
York.

\bibitem{vc}
Vapnik, V. N. and Chervonenkis, A. Ja. (1974). {\it Theory of
Pattern Recognition}. Nauka, Moscow (in Russian).

\bibitem{Yan99}
Yang, Y. (1999) Minimax nonparametric classification - Part I: Rates
of convergence, Part II: Model selection for adaptation. {\it IEEE
Trans. Inf. Theory} {\bf 45} 2271-2292.


\end{thebibliography}
\end{document}